\documentclass[11pt,twoside]{amsart}

\usepackage{amsfonts}
\usepackage{amssymb}
\usepackage{amscd}
\usepackage[dvips]{graphicx}
\usepackage{xcolor}

\usepackage[all]{xy}

\title[DG-orders: class groups, dg-radicals, maximal dg-orders]{Differential graded
orders, their class groups and id\`eles}

\date{January 9, 2023; revised September 10, 2024}

\author{Alexander Zimmermann}
\address{
Universit\'e de Picardie,
D\'epartement de Math\'ematiques et LAMFA (UMR 7352 du CNRS),
33 rue St Leu,
F-80039 Amiens Cedex 1,
France}
\email{alexander.zimmermann@u-picardie.fr}

\newtheorem{Lemma1}{{Lemma}}[section]
\newtheorem{Def1}[Lemma1]{{Definition}}
\newtheorem{Prop1}[Lemma1]{{Proposition}}
\newtheorem{Claim1}[Lemma1]{{Claim}}
\newtheorem{Rem1}[Lemma1]{{Remark}}
\newtheorem{Cor1}[Lemma1]{{Corollary}}
\newtheorem{Ex1}[Lemma1]{{Example}}
\newtheorem{Qu1}[Lemma1]{{Question}}
\newtheorem{Theo1}[Lemma1]{{Theorem}}
\newtheorem*{Theo2}{{Theorem}}

\newenvironment{Lemma}{\begin{Lemma1}}{\end{Lemma1}}
\newenvironment{Def}{\begin{Def1}\em}{\end{Def1}}
\newenvironment{Prop}{\begin{Prop1}}{\end{Prop1}}
\newenvironment{Rem}{\begin{Rem1}\rm}{\end{Rem1}}
\newenvironment{Theorem}{\begin{Theo1}}{\end{Theo1}}

\newenvironment{Cor}{\begin{Cor1}}{\end{Cor1}}

\newenvironment{Example}{\begin{Ex1}\em}{\end{Ex1}}
\newenvironment{Question}{\begin{Qu1}\em}{\end{Qu1}}

\newcommand{\lra}{\longrightarrow}

\newcommand{\ra}{\rightarrow}
\newcommand{\sdp}{\times\kern-.2em\vrule height1.1ex depth-.05ex}
\newcommand{\epi}{\lra \kern-.8em\ra}

\newcommand{\Q}{{\mathbb Q}}
\newcommand{\N}{{\mathbb N}}

\newcommand{\ol}{\overline}
\newcommand{\Sn}{{\mathfrak S}}

\newcommand{\Z}{{\mathbb Z}}

\newcommand{\dickebox}{{\vrule height5pt width5pt depth0pt}}

\newcommand{\im}{\textup{im}}
\newcommand{\dgrad}{\textup{dgrad}}
\newcommand{\ann}{\textup{ann}}

\newcommand{\red}{}

 \normalbaselines
\begin{document}

\begin{abstract}
For a Dedekind domain $R$ with field of fractions $K$ a classical $R$-order in a semisimple $K$-algebra
$A$ is an $R$-projective $R$-subalgebra $\Lambda$ of $A$ such that $K\Lambda=A$. We study
differential graded $K$-algebras which are semisimple as $K$-algebras and define
differential graded $R$-orders as a differential graded $R$-subalgebras, which are
in addition classical $R$-orders in $A$. We give a series of examples for such differential graded
algebras and orders. We show that any differential graded
$R$-order is contained in a maximal differential graded order.
We develop parts of the classical
ring theory in the differential graded setting, in particular the properties of
analogues of the Jacobson radical. We further define class groups of differential graded orders
as subgroups of the Grothendieck group of locally free differential graded modules. We define
id\`eles in this setting showing that these id\`ele groups {\red map surjectively} to the differential graded
class group. Finally we give a homomorphism to the class group of the homology of the differential
graded order and prove a Mayer-Vietoris like sequence for each central idempotent of $A$, including
the analogous one for the kernel groups of these morphisms.
\end{abstract}

\maketitle

\section*{Introduction}

Differential graded algebras were introduced by Cartan in \cite{Cartandg}. These are, by definition,
a pair $(A,d)$ of an associative $\Z$-graded algebra $A$ with an endomorphism $d:A\lra A$
of degree $1$ and $d^2=0$, satisfying $d(a\cdot b)=a\cdot d(b)+(-1)^{|a|}a\cdot d(b)$
for all homogeneous elements $a,b\in A$. Differential graded algebras became prominent in the past
few decades in order to provide a powerful tool for homological algebra.

Most research was directed in this perspective. In particular, as differential graded
algebras occur naturally in the homological algebra of algebraic and differential
varieties, these examples provided the guideline of reasonable assumptions. In particular,
a quite systematic assumption is that differential graded algebras are connective, that is
all homogeneous components in positve degrees are assumed to vanish. Often, they
are algebras of infinite dimension over a base field.

We ask what happens if the algebras are more linked to classical finite dimensional
algebras. Differential graded algebras are algebras at first, and so one might be
tempted to understand
ring theoretic properties of differential graded algebras. To our
greatest surprise we could not find much research which was carried out up to now in this direction.
We fist study properties such as Jacobson radicals or analogues of these, and
see that the concepts in this direction in the differential graded case are different.
Unlike the classical situation the intersection of all
differential graded maximal left ideals is not the same as the intersection
of annihilators of differential graded simple left modules. This holds however up to homotopy.
Further, we prove a Nakayama's lemma
for the two-sided notion given by intersection of annihilators of simple dg-modules.

We then study differential graded orders in differential graded semisimple algebras.
Since {\red in characteristic different from $2$}
the differential vanishes automatically on central idempotents, we study
finite dimensional simple algebras, i.e. matrix algebras over skew fields by
Wedderburn's theorem. We provide examples that there are various interesting differential graded
structures on these matrix algebras. They bear interesting properties.
In particular, their homology is not semisimple in general.

Then, we consider differential graded orders in these algebras, naively as classical orders which are
stable under the differential graded structure. Graded orders were studied by LeBruyn, van Oystaen,
and van den Bergh in \cite{GradedOrders}. The differential gives an important extra
structure. One may ask if it would not be reasonable to ask in addition that the
homology is an order in the homology of the algebra. However, this additional condition is
restrictive and it does not seem to allow a rich theory.

We prove that maximal differential graded orders always exist, that any differential order
is contained in a maximal one, and give examples that not all maximal orders allow a
differential structure. Clearly, our
algebras are not connective. Non connective differential graded algebras
may have strange properties. We mention Raedschelder and Stevenson \cite{Raedschelder} for
striking examples and properties in this direction.

Finally we define class groups of locally free differential graded orders. We give a
$K$-theoretic definition and prove its equivalence with a second definition provided
by id\`eles. The only modification which is necessary, is that we need to consider
the subalgebra of cycles, and not the entire algebra. We then show that the
functor 'taking homology' then maps to the ordinary class group of the homology algebra,
at least if the homology algebra of the ambiant semisimple algebra is semisimple.
The kernel of the resulting map then bounds the fibre of stable quasi-isomorphism classes
instead of stable isomorphism classes of locally free modules.

We finally mention that a differential graded Jordan-Zassenhaus theorem is not available
for the moment. Hence I do not known if  differential graded class groups are
finite.

Here is an outline of the paper. In Section~\ref{generalities} we review some basic
facts on differential graded algebras and differential graded  modules.
{\red In Section~\ref{generalitiesgraded} we recall a result due to Aldrich and Garcia Rozas
on semisimplicity of the category of differential graded modules.}
In Section~\ref{Examplesection} we provide basically two classes of algebras which
are simple as algebras and differential graded. Choosing special cases provide a
wealth of interesting examples for what follows. Section~\ref{ringtheoryofdgalgebras}
develops classical ring theory of differential graded algebras, in particular
different concepts of radicals of differential graded algebras, providing the technical tools
for discussions in further sections. Section~\ref{diffgradedordersection} then introduces
our main object of study, differential graded orders and differential graded lattices.
We provide a weak version and a strong
version in parallel, but focus mainly on the weak version, which seems to have
more and easier accessible properties.
As fist property we study maximal differential graded orders in Section~\ref{maximaldiffgradedorders}.
In Section~\ref{classgroupsection} we introduce class groups of differential graded orders
and differential graded id\`eles, show the equivalence of the K-theoretic and the id\`ele
theoretic setting.  In Section~\ref{mayervietorissection} we prove a Mayer-Vietoris like
sequence for class groups of differential graded orders.

\section{Foundations of dg-algebras and dg-modules}

\label{generalities}

\subsection{Generalities}

Let $R$ be a commutative ring. Recall from Cartan \cite{Cartandg} and Keller \cite{Kellerdgandtilting} that
\begin{itemize}
\item
a differential graded $R$-algebra (or dg-algebra for short) is a $\Z$-graded
$R$-algebra $A=\bigoplus_{n\in\Z}A_n$ together with a graded $R$-linear
endomorphism $d$ of square $0$
of degree $1$ (i.e. $d(A_n)\subseteq A_{n+1}$ and $d\circ d=0$) and such that $$d(ab)=d(a)b+(-1)^{|a|}ad(b)$$
for all homogeneous elements  $a,b\in A$. A homomorphism of differential graded algebras
$f:(A,d_A)\ra (B,d_B)$ is a degree $0$ homogeneous $R$-algebra map $f$ such that
$f\circ d_A=d_B\circ f$.
{\red\item
If $(A,d)$ is a differential graded algebra, then $(A^{op},d^{op})$ is a differential
graded algebra (cf e.g. \cite[Definition 11.1]{Stacks})
with $x\cdot_{op}y:=(-1)^{|x|\cdot|y|}yx$ for any homogeneous elements $x,y\in A$, and $d^{op}(x)=d(x)$.
We hence write $d^{op}=d$.}
\item
A differential graded left $A$-module (or dg-module for short) is then an $A$-module $M$, $\Z$-graded as an
$R$-module with graded $R$-linear
endomorphism $d_M$ of square $0$ and of degree $1$, such that $$d_M(ma)=d_M(m)a+(-1)^{|m|}md(a)$$
for all homogeneous elements $a\in A$ and $m\in M$.
{\red A differential graded $(A,d)$-right module is a differential graded $(A^{op},d)$-left module.}
\item
Let $(A,d_A)$ be a differential graded $R$-algebra and let $(M,\delta_M)$ and $(N,\delta_N)$ be
differential graded $(A,d_A)$-modules. Then a homomorphism of differential
graded modules is an $R$-linear map $f:M\rightarrow N$, homogeneous of degree $0$
with $f\circ \delta_M=\delta_N\circ f$, with
$f(am)=af(m)$ for all $a\in A$ and $m\in M$.
\item
If $(M,d_M)$ is a dg-$A$-module, then let $M[k]$ be the dg-module given by $(M[k])_n=M_{k+n}$
for all $n\in\Z$ and $d_{M[k]}=d_M$. It is easy to see that $M[k]$ is a dg-module again.
\item Let $(M,d_M)$ and $(N,d_N)$ be
differential graded $(A,d)$-modules.
The homomorphism complex $Hom_A^\bullet(M,N)$ is the $\Z$-graded $A$-module given by
{\red\begin{align*}
(Hom_A^\bullet(M,N))_n:=\{f:M\ra N\;|&\;f\in Hom_\Z(M,N)\;\mbox{ and }\\
&f(M_k)\subseteq N_{k+n}\textup{ and }\\
&f(am)=(-1)^{|a|n}a\cdot f(m)\}.
\end{align*}}
The elements $f$ of $Hom_A^\bullet(M,N)$  are not asked to be compatible with the differentials in any way.
Let $d_{Hom}: Hom_A^\bullet(M,N)\ra Hom_A^\bullet(M,N)$ given by
$$d_{Hom}(f):=d_N\circ f-(-1)^{|f|}f\circ d_M.$$
Then
\begin{eqnarray*}d_{Hom}^2(f)&=&d_{Hom}(d_N\circ f-(-1)^{|f|}f\circ d_M)\\
&=&
d_N\circ(d_N\circ f-(-1)^{|f|}f\circ d_M)\\
&&-(-1)^{|f|-1}(d_N\circ f-(-1)^{|f|}f\circ d_M)\circ d_M\\
&=&
(-1)^{{\red |f|+1}}(d_N\circ f\circ d_M-d_N\circ f\circ d_M)\\
&=&0
\end{eqnarray*}
Therefore $(Hom_A^\bullet(M,N),d_{Hom})$ is a complex of $R$-modules.
\item The category $A-dgMod$ of dg-modules over the dg-algebra $(A,d)$ is abelian with morphisms the
degree zero cycles of $(Hom_A^\bullet(M,N),d_{Hom})$ (cf~\cite[Lemma 22.4.2]{Stacks}).
If $A$ is Noetherian, then the subcategory
$A-dgmod$ of finitely generated dg $A$ modules is abelian as well.
\item The homotopy category $K(A-dgMod)$ is the category with objects $A$ dg-modules and
morphisms $Hom_{K(A-dgMod)}(M,N)=H_0(Hom_A^\bullet(M,N),d_{Hom})$.
\item If $A$ is artinian and $(A,d)$ a differential graded algebra. Then the category
$A-dgmod$ has split idempotents and is hence Krull-Schmidt. This follows from the fact that
the image of a homomorphism of dg-modules is a dg-module again.
\end{itemize}

\begin{Lemma}\label{Hombullet} Let $(A,d)$ be a differential graded $R$-algebra, and let
$(M,d_M)$ and $(N,d_N)$ be differential graded $(A,d)$-modules.
\begin{itemize}
\item
Then $Z_0(Hom_A^\bullet(M,N),d_{Hom})=({\ker \red (d_{Hom})_0}$ consists the homomorphisms
of complexes, and $H_0(Hom_A^\bullet(M,N),d_{Hom})$
is the morphisms $M\ra N$ in the homotopy category.
\item If $M=N$, then $End^\bullet_A(M)$ is a
dg-algebra and composition of graded morphisms yields a differential
graded $End^\bullet_A(M)$-module structure on $Hom^\bullet_A(M,N)$.
\end{itemize}
\end{Lemma}

Proof: \cite{Kellerdgandtilting} \dickebox

\begin{Lemma} \label{oneisinthecycles}
$d_A(1)=0$.
\end{Lemma}

Proof.
$d_A(1)=1\cdot d_A(1)+d_A(1)\cdot 1$, which implies
$d_A(1)=0$.
%
\dickebox

\begin{Cor} \label{cyclesaresubalgebras} Let {\red $R$ be a commutative ring} and
let $(A,d)$ be a differential graded $R$-algebra. Then $\ker(d)$ is a graded subalgebra of $A$.
\end{Cor}

Proof. Since $d$ is $R$-linear, $\ker(d)$ is an $R$-module.
Since $d$ is homogeneous, $\ker(d)$ is graded. By Lemma~\ref{oneisinthecycles} $\ker(d)$ contains $1$.
Further, $d(xy)=d(x)y\pm xd(y)$ for any homogeneous elements $x,y$ shows that
$x,y\in\ker(d)\Rightarrow xy\in\ker(d).$
\dickebox

\begin{Rem}
We note that Example~\ref{Examplediffgradedalgebra} shows that in general $d(u)\neq 0$ for a unit $u$ in $A$.
Further, units are not necessarily in degree $0$, as is shown by Example~\ref{Examplediffgradedalgebra}.
\end{Rem}

\begin{Lemma}\label{dgalgebrainducesalgebraonH}
Let $(A,d)$ be a differential graded algebra. Then, the algebra structure of $A$ induces a
graded algebra structure on $H(A,d)$ (and hence differential graded with differential $0$).
\end{Lemma}

Proof. Let $a,b\in\ker(d)$ {\red be homogeneous.} Then {\red $d(a\cdot b)=d(a)\cdot b+(-1)^{|a|} a\cdot d(b)=0.$}
Hence $a\cdot b\in\ker(d)$. Let $a=d(x)$ and $c\in \ker(d)$ be homogeneous elements.
Then {\red $$ac=d(x)c= d(xc)-(-1)^{|b|} xd(c)= d(xc)\in \im(d).$$}
Likewise $ca\in\im(d)$. Hence,
the multiplication is well-defined. The fact that the additive law is
well-defined and combines to a ring structure is trivial. \dickebox

\begin{Lemma}\label{homologymodule}
Let $(A,d)$ be a differential graded algebra and let $(M,d_M)$ be a
differential graded {\red right (resp. left)} $A$-module. Then the $A$-module structure on $M$
induces a {\red right (resp. left)}  $H(A,d)$-module structure on $H(M,d_M)$.
\end{Lemma}

Proof. We only need to show that the induced action is well-defined. {\red We shall
give the proof for a right module. The left module case is analogous. }

Let $a\in\ker(d)$ be a homogeneous element, and let $m\in\ker(d_M)$
be a homogeneous element. Then
$$d_M(m\cdot a)=d_M(m)\cdot a+(-1)^{|m|}m\cdot d(a)=0$$
Hence $m\cdot a\in\ker(d_M)$ again.
Moreover,
\begin{eqnarray*}
d_M(n)\cdot(a+d(b))&=&d_M(n)a+d_M(n)d(b)\\
&=&
d_M(na)-(-1)^{|a|}nd(a)+d_M(nd(b))-(-1)^{|n|}nd^2(b)\\
&=&d_M(na)+d_M(nd(b))-(-1)^{|a|}nd(a)\\
&=&d_M(na)+d_M(nd(b))
\end{eqnarray*}
is in $\im(d_M)$, since $a\in\ker(d)$ and hence $nd(a)=0$.
Therefore, the multiplication of $H(A)$ on $H(M)$ is well-defined.
\dickebox

\section{Semisimplicity in the category of differential graded modules}

\label{generalitiesgraded}

We recall a theorem of Tempest and Garcia-{\red Rozas}~\cite{Tempest-Garcia-Rochas}. A differential graded
$A$-module $X$ is simple if the only differential graded $A$-submodules are $0$ and $X$.
{\red An abelian category is semisimple if all subobjects admit a complement.}

\begin{Theorem} \label{GarciaRozastheorem}
\begin{itemize}
\item \cite[Theorem 4.7]{Tempest-Garcia-Rochas}
Let $(A,d)$ be a differential graded algebra. Then the following are equivalent:
\begin{enumerate}
\item The regular differential graded $A$-module
$A$ is projective in the category of dg-modules
\item $A$ is {\red acyclic}. 
\item $1_A\in im(d)$
\item any left dg-module is {\red acyclic}. 
\item The functor $A\otimes_{Z(A)}-$ from  graded $Z(A)$-modules to differential graded $A$-modules is an equivalence with quasi-inverse $Z(-)$.
\end{enumerate}
\item \cite[Theorem 5.3]{Tempest-Garcia-Rochas}\label{tempest-GarciaRozas}
{\red The category of differential graded modules over $(A,d)$ is semisimple precisely when
the category of graded modules over the graded algebra of cycles $Z(A)=\ker(d)$ is
semisimple and $A$ is an {\red acyclic } 
complex. This in turn is equivalent with the left regular differential graded $(A,d)$-module $(A,d)$ is
semisimple as a differential graded module, in the sense that all differential graded ideals
have a complement}.
\end{itemize}
\end{Theorem}

Note that Tempest and Garcia-Rozas show \cite[Lemma 4.2]{Tempest-Garcia-Rochas} that
if $1_A=d(z)$, then $A=Z(A)\oplus Z(A)y$.

\bigskip

\begin{Rem}\label{blocksofdgalgebras}
Let $K$ be a commutative ring {\red such that $2$ is a regular element in $K$,} and let
$(A,d)$ be a differential graded $K$-algebra with centre $Z(A)$,
and let $e^2=e\in Z(A)$ be homogeneous.
Then, {\red copying the argument of \cite[Theorem 1.4 and Corollary 1.5]{DascalescuIon}, since $e^2=e\in Z(A)$, }
we get that $e$ is necessarily in $A_0$. Moreover,
$$d(e)=d(e^2)=d(e)e+ed(e)=2ed(e)=4ed(e)
$$
Hence $ed(e)=0$. Moreover $(1-e)d(e)=(1-e)d(e)e+(1-e)ed(e)=0$, and hence $d(e)=0$.
Therefore, if $A=A_1\times A_2$ as algebra,
then $(A,d)=(A_1,d)\times (A_2,d)$ as differential graded algebras.
\end{Rem}

\begin{Cor}
As we have seen in Remark~\ref{blocksofdgalgebras} for a differential
graded algebra $A$ which is semisimple artinian as algebra {\red over a field of
characteristic different from $2$}, by Wedderburn's theorem
we have that $A$ is a direct sum of differential graded algebras, each of which is a
matrix algebra of a skew field.
\end{Cor}

\section{Some examples}

We first give an example for a $\Z$-grading of matrix algebras.

\begin{Example} \label{gradingsofmatrixrings}
\begin{enumerate}
\item \label{gradingsofmatrixrings-1}
Let $A=Mat_{n\times n}(D)$ for some skew-field $D$, finite dimensional over $K$, and $n>1$
an integer. Then $A$ is graded by putting $A_0$ the subalgebra given by the main
diagonal entries
$$\left(\begin{array}{ccccc}\ast&0&\dots&\dots&0\\0&\ast&0&&\vdots\\
\vdots&\ddots&\ast&\ddots&\vdots\\
\vdots&&\ddots&\ddots&0\\
0&\dots&\dots&0&\ast\end{array}\right)$$
Then upper one diagonal coefficients $(a_{i,i+1})_{i}$
give the degree $1$
$$\left(\begin{array}{ccccc}0&\ast&0&\dots&0\\0&0&\ast&\ddots&\vdots\\
\vdots&\ddots&\ddots&\ddots&0\\
\vdots&&\ddots&0&\ast\\
0&\dots&\dots&0&0\end{array}\right)$$ the upper $2$ diagonal coefficients $(a_{i,i+2})_{i}$  give the
degree $2$, etc. The lower one diagonal coefficients $(a_{i+1,i})_{i}$
give the degree $-1$, the lower $2$ diagonal coefficients $(a_{i+2,i})_{i}$
give the degree $-2$, etc. Elementary matrix multiplication proves the definition for a
graded algebra. Moreover, $A$ is semisimple as an algebra.
\item
Let $A=\bigoplus_{n\in\Z}A_n$ be a graded semisimple algebra, then for any $u\in A^\times$
we see that $A^u:=\bigoplus_{n\in\Z}uA_nu^{-1}$ is again a graded semisimple algebra.
\end{enumerate}
\end{Example}

\label{Examplesection}

We produce a series of examples showing that our concept is non trivial and allows
interesting phenomena.

\begin{Rem} \label{AenotinAde}
Note that Example~\ref{gradingsofmatrixrings} provides gradings such that
$A\cdot e\not\subseteq A\cdot d(e)$ for any primitive idempotent $e$ of degree $0$
and any differential $d$ of degree $1$. Also for
$End_K^\bullet(L)$ for a bounded complex $L$ of $K$-vector spaces there is such an idempotent.
Indeed, in this case the degree $0$
component is formed by matrix algebras along the diagonal, and the degree $1$ component is
formed by matrices right up to these.
\end{Rem}

\begin{Prop}\label{progenerator}
Let $K$ be a field and let $(A,d)$ be a finite dimensional differential graded
algebra. Suppose that $A$ is a split simple $K$-algebra.
Suppose that there is a primitive idempotent $e$ of $A$ such that $A\cdot e\not\subseteq A\cdot d(e)$.
Then there is a  bounded
complex $L$ of $K$-modules such that $A\simeq Hom_K^\bullet(L,L)$ as differential graded algebras.
Conversely,  $A=Hom_K^\bullet(L,L)$ is differential graded, finite dimensional simple as algebra, such that there is a primitive idempotent $e$ with
$A\cdot e\not\subseteq A\cdot d(e)$.
\end{Prop}

Proof. If $L$ is a bounded complex of $K$-vector spaces, then
$Hom_K^\bullet(L,L)$ is a full matrix ring over $K$, as ungraded algebra, and hence simple
as algebra. Further, $Hom_K^\bullet(L,L)$ is a differential graded algebra by
Lemma~\ref{Hombullet}.

Conversely, let $K$ be a field and let $(A,d)$ be a finite dimensional differential graded
algebra. Suppose that $A$ is a split simple $K$-algebra. By Wedderburn's theorem, $A$ is a
full  matrix algebra over $K$. Let $e$ be a primitive idempotent of $A$ satisfying
$Ae\not\subseteq Ad(e)$. Then
$$M:=A\cdot e+A\cdot d(e)\textup{ and }N:=A\cdot d(e)$$
are differential graded $(A,d)$-modules. Further $N< M$ and $$L:=M/N\neq 0$$ is a
differential graded $(A,d)$-module. As $A$-module, we see that $L\simeq Ae$ is a progenerator.
Hence $L$ is a natural differential graded $(A,d)-(End^\bullet(L),d_{Hom})$
bimodule. Now, for any homogeneous $a\in A$, left multiplication by $a$ gives a
homogeneous element $\varphi(a)\in End^\bullet(L)$. Further, $\varphi$ is additive, sends
$1\in A$ to the identity on $L$, and induces a ring homomorphism
$$\varphi:A\lra  End^\bullet(L).$$
Since $L$ is a progenerator, $\varphi$ is injective. Since
$\dim_K(A)=\dim_K(End^\bullet(L))$, we get that $\varphi$ is an isomorphism of algebras.
Now, for any homogeneous $a,b\in A$, we have
$$d(a)b=d(ab)-(-1)^{|a|}ad(b)$$
we get
$$\varphi(d(a))=d\circ\varphi(a)-(-1)^{|a|}\varphi(a)\circ d=d_{Hom}(\varphi(a))$$
and therefore $\varphi$ is an isomorphism of differential graded algebras.
Further, by Remark~\ref{AenotinAde} there is a primitive idempotent $e$ of degree $0$,
by degree considerations we get $Ae\not\subseteq Ad(e)$.
\dickebox

{\red \begin{Rem}\label{strangehypothesisisunnecessary}
The hypothesis that there is a primitive idempotent $e$ of $A$ such that
$A\cdot e\not\subseteq A\cdot d(e)$ is superfluous. Indeed,
\cite[Corollary 1.5]{DascalescuIon} shows that $A$ is isomorphic,
as a graded algebra, to a matrix algebra with the main diagonal in degree $0$.
Transporting the differential structure via this isomorphism, we obtain the existence
of such an idempotent.
\end{Rem}}

{\red\begin{Rem}
After having finished and submitted the manuscript I discovered that D. Orlov
defined our notion of simple differential graded algebra earlier in \cite{Orlov1},
and called it abstractly simple.
Moreover, he proved the statement of Proposition~\ref{progenerator} by completely different means.
His proof uses scheme theoretic arguments. However, he has to assume that the primitive
central idempotents are in degree $0$. Our approach however gives that this
can be assumed to be automatically satisfied
using \cite[Corollary 1.5]{DascalescuIon}.
\end{Rem}}

\begin{Example}\label{Examplediffgradedalgebra}
We consider a field $K$ and the algebra of $2\times 2$ matrices over $K$.
We shall use Proposition~\ref{progenerator}.
Except the stalk complex and complexes with differential $0$, {\red up to shift},
the only possible complex realising this dg-algebra is
$$\cdots\lra 0\lra K\stackrel{\cdot x}{\lra}K\lra 0\lra\cdots,$$
{\red concentrated in degree $0$ and $1$. We denote by $\delta$ the
differential on the complex. Endomorphisms of degree $0$ are given by two scalars $u$ and $v$
$$
\xymatrix{
K\ar[r]^{\cdot x}\ar@(l,u)[]^{\cdot u}&K\ar@(r,u)[]_{\cdot v}.
}
$$
A morphism of degree $1$ will map the degree $0$ homogeneous component to the degree $1$ homogeneous component, i.e. is given by multiplication by a scalar $w$.
A morphism of degree $-1$ will map the degree $1$ homogeneous component to the degree $0$ homogeneous component, i.e. is given by multiplication by a scalar $z$.
$$
\xymatrix{\textup{degree $1$ morphism}&&&&\textup{degree $-1$ morphism}\\
K\ar[r]^{\cdot x}\ar[rd]^{\cdot w}&K&&K\ar[r]^{\cdot x}&K\ar[ld]_{\cdot z}\\
K\ar[r]^{\cdot x}&K&&K\ar[r]^{\cdot x}&K}
$$
Then, these morphism correspond to a matrix
$$\left(\begin{array}{cc}u&w\\z&v\end{array}\right)\in Mat_{2\times 2}(K)$$
with the names of the variables chosen as in the above maps indicating the identification
of the endomorphisms with the matrices, and the grading given as in
Example~\ref{gradingsofmatrixrings}.\ref{gradingsofmatrixrings-1}.
Here, in order to respect the usual multiplication of matrices and
endomorphisms of the complex, we need to write the maps on the right, and compose them accordingly.

Compute the differential.
Consider first a degree $-1$ morphism $\alpha_z$. Then, using that the maps apply on the right and hence
$d(\gamma)=\gamma d-(-1)^{|\gamma|}d\gamma$ as mappings acting on the complex. Hence
$$
d(\alpha_z)=\alpha_z\delta-(-1)^{|\alpha_z|}\delta\alpha_z=
\left(\begin{array}{cc}xz&0\\0&0\end{array}\right)+\left(\begin{array}{cc}0&0\\0&xz\end{array}\right)=x\cdot\textup{id}.
$$
Then, for the differential of a degree $0$ morphism $\beta_{(u,v)}$ acting as $u$ on the degree $0$ component and
as $v$ on the degree $1$-component we get
\begin{eqnarray*}
d(\beta_{u,v})&=&\beta_{u,v}\delta-(-1)^{|\beta_{u,v}|}\delta\beta_{u,v}\\
&=&\beta_{u,v}\delta-\delta\beta_{u,v}\\
&=&
\left(\begin{array}{cc}0&-xu\\0&0\end{array}\right)+\left(\begin{array}{cc}0&xv\\0&0\end{array}\right)\\
&=&
\left(\begin{array}{cc}0&x(v-u)\\0&0\end{array}\right)
\end{eqnarray*}
Note that this result is a consequence of the identities
$$\left(\begin{array}{cc}0&0\\1&0\end{array}\right)\cdot  \left(\begin{array}{cc}0&1\\0&0\end{array}\right)
=\left(\begin{array}{cc}0&0\\0&1\end{array}\right)\textup{ and }
  \left(\begin{array}{cc}0&1\\0&0\end{array}\right)\cdot \left(\begin{array}{cc}0&0\\1&0\end{array}\right)
=\left(\begin{array}{cc}1&0\\0&0\end{array}\right)
$$
together with Leibniz formula and the result of $d(\alpha_z)$.
}
This gives a dg-algebra and any differential is of this form for
some $x\in K$. If $x\neq 0$, the homology is $0$.

Note that if $x\neq 0$, then the kernel of the differential is
$${\red \ker(d)}=\{\left(\begin{array}{cc}a&b\\0&a\end{array}\right)\;|\;a,b\in K\}\simeq K[\epsilon]/\epsilon^2$$
with $\epsilon$ in degree $1$.
This is not semisimple, and hence we do not get a semisimple differential graded algebra in the sense of
Aldrich and Garcia-{\red Rozas}~\cite{Tempest-Garcia-Rochas}. However, for $x$ invertible we have $z=\frac{1}{x}\left(\begin{array}{cc}0&0\\1&0\end{array}\right)$ is a preimage of $1$, as
required in Aldrich and Garcia-{\red Rozas}~\cite{Tempest-Garcia-Rochas} structure theorem.

We further observe that $\left(\begin{array}{cc}0&1\\1&0\end{array}\right)$ is invertible, non homogeneous
with a summand of degree $-1$ and another summand of degree $1$.
However,
$$
d(\left(\begin{array}{cc}0&1\\1&0\end{array}\right))=d(\left(\begin{array}{cc}0&1\\0&0\end{array}\right))+
d(\left(\begin{array}{cc}0&0\\1&0\end{array}\right))=\left(\begin{array}{cc}x&0\\0&x\end{array}\right)
$$
is non zero for $x\neq 0$. Hence, the differential of invertible elements are not necessarily $0$.
Nevertheless, $d(1)=0$.

As any differential graded $(A,d)$-module $(M,\delta)$ is at first an $A$-module, it is useful to
consider the possible {\red left} dg-module structures on $K^2$ for the above dg-algebra $(Mat_2(K),d_x)$. So, let
$\delta(\left(\begin{array}{c}0\\1\end{array}\right))=\left(\begin{array}{c}u\\v\end{array}\right)$ for some $u,v\in K$ and $\delta(\left(\begin{array}{c}1\\0\end{array}\right))=\left(\begin{array}{c}s\\t\end{array}\right)$ for some $s,t\in K$.
\begin{eqnarray*}
0&=&\delta(\left(\begin{array}{cc}0&0\\1&0\end{array}\right)\cdot \left(\begin{array}{c}0\\1\end{array}\right))
=\left(\begin{array}{cc}x&0\\0&x\end{array}\right)\cdot \left(\begin{array}{c}0\\1\end{array}\right)-
\left(\begin{array}{cc}0&0\\1&0\end{array}\right)\cdot \left(\begin{array}{c}u\\v\end{array}\right)
=\left(\begin{array}{c}0\\ x-u\end{array}\right)\\
0&=&\delta(\left(\begin{array}{cc}1&0\\0&0\end{array}\right)\cdot \left(\begin{array}{c}0\\1\end{array}\right))
=\left(\begin{array}{cc}0&-x\\0&0\end{array}\right)\cdot \left(\begin{array}{c}0\\1\end{array}\right)+
\left(\begin{array}{cc}1&0\\0&0\end{array}\right)\cdot \left(\begin{array}{c}u\\v\end{array}\right)
=\left(\begin{array}{c}u-x\\ 0\end{array}\right)\\
\left(\begin{array}{c}u\\ v\end{array}\right)&=&\delta(\left(\begin{array}{cc}0&0\\0&1\end{array}\right)\cdot \left(\begin{array}{c}0\\1\end{array}\right))
=\left(\begin{array}{cc}0&x\\0&0\end{array}\right)\cdot \left(\begin{array}{c}0\\1\end{array}\right)+
\left(\begin{array}{cc}0&0\\0&1\end{array}\right)\cdot \left(\begin{array}{c}u\\v\end{array}\right)
=\left(\begin{array}{c}x\\ v\end{array}\right)\\
\left(\begin{array}{c}s\\ t\end{array}\right)&=&\delta(\left(\begin{array}{cc}0&1\\0&0\end{array}\right)\cdot \left(\begin{array}{c}0\\1\end{array}\right))
=
-\left(\begin{array}{cc}0&1\\0&0\end{array}\right)\cdot \left(\begin{array}{c}u\\v\end{array}\right)
=\left(\begin{array}{c}-v\\ 0\end{array}\right)\\
\left(\begin{array}{c}u\\ v\end{array}\right)&=&\delta(\left(\begin{array}{cc}0&0\\1&0\end{array}\right)\cdot \left(\begin{array}{c}1\\0\end{array}\right))
=
\left(\begin{array}{cc}x&0\\0&x\end{array}\right)\cdot \left(\begin{array}{c}1\\0\end{array}\right)-
\left(\begin{array}{cc}0&0\\1&0\end{array}\right)\cdot \left(\begin{array}{c}s\\t\end{array}\right)
=\left(\begin{array}{c}x\\ -s\end{array}\right)\\
\left(\begin{array}{c}as\\ 0\end{array}\right)&=&\delta(\left(\begin{array}{cc}a&0\\0&b\end{array}\right)\cdot \left(\begin{array}{c}1\\0\end{array}\right))
=
\left(\begin{array}{cc}0&x(b-a)\\0&0\end{array}\right)\cdot \left(\begin{array}{c}1\\0\end{array}\right)+
\left(\begin{array}{cc}a&0\\0&b\end{array}\right)\cdot \left(\begin{array}{c}s\\t\end{array}\right)
=\left(\begin{array}{c}as\\ bt\end{array}\right)\\
0&=&\delta(\left(\begin{array}{cc}0&1\\0&0\end{array}\right)\cdot \left(\begin{array}{c}1\\0\end{array}\right))
=
-\left(\begin{array}{cc}0&1\\0&0\end{array}\right)\cdot \left(\begin{array}{c}s\\t\end{array}\right)
=\left(\begin{array}{c}-t\\ 0\end{array}\right)\\
0&=&\delta^2(\left(\begin{array}c 0\\1\end{array}\right))=\delta(\left(\begin{array}c x\\v\end{array}\right))
=\left(\begin{array}c vx\\v^2\end{array}\right)+\left(\begin{array}c vx\\0\end{array}\right)=
\left(\begin{array}c 2vx\\v^2\end{array}\right)
\end{eqnarray*}
Hence  $s=v=t=0$,  $x=u$. Since $\delta$ is of degree $1$, we need to have
$\left(\begin{array}{c}1\\0\end{array}\right)$ is of degree $n$ and
$\left(\begin{array}{c}0\\1\end{array}\right)$ is of degree $n-1$ for some integer $n$.
Up to shift there is a unique dg-module structure given by $\delta(\left(\begin{array}c p\\q\end{array}\right))=\left(\begin{array}{c}xq\\0\end{array}\right)$.
Of course, this is a simple differential graded module.

Let $(C^*,\widehat\delta)$ be a bounded complex of finitely generated
$Mat_{2\times 2}(K)$-modules. Each homogeneous component
$C^n$ is isomorphic to a a direct sum of differential graded modules $(K^2,\delta)$ as above.
Then, the total complex is a differential graded module. Using the above explicit computation
can be used to show that this gives again the only possible differential graded structures on the
corresponding graded module.
\end{Example}

\begin{Rem}
I wish to thank Bernhard Keller who noted\footnote{email from December 24, 2022 to the author}
that the above differential graded algebra is the differential graded
endomorphism algebra $Hom_K^\bullet(M,M)$ for $M$ being the complex
$$\cdots \lra 0\lra K\stackrel{\cdot x}\lra K\lra 0\lra \cdots.$$
\end{Rem}

\begin{Example} \label{dgalgebraonthreebythree}
We consider the $3\times 3$ matrix ring over a commutative ring $R$. We will use an alternative grading coming from the fact that this ring is Morita equivalent to the $2\time 2$ matrix ring over $R$.
\begin{eqnarray*}
Mat_3(R)=\{\left(\begin{array}{ccc}
a_{11}&a_{12}&a_{13}\\ a_{21}&a_{22}&a_{23}\\ a_{31}&a_{32}&a_{33}
\end{array}\right)\;&|&\; \red{deg(a_{31})=deg(a_{32})=}-1; \\
&&deg(a_{11})=deg(a_{12})=deg(a_{21})=deg(a_{22})=deg(a_{33})=0;\\
&&deg(a_{13})=deg(a_{23})=1\}
\end{eqnarray*}
This grading comes from the complex $L$
$$
\cdots \lra 0\lra R\stackrel{a_{11}\choose a_{21}}\lra R^2\lra 0\lra \cdots
$$
We make explicit the differentials.

We get
$$
d(\left(\begin{array}{ccc}
0&0&0\\ 0&0&0\\ x&y&0
\end{array}\right))=
\left(\begin{array}{ccc}
a_{11}x&a_{11}y&0\\ a_{21}x&a_{21}y&0\\ 0&0&a_{11}x+a_{21}y
\end{array}\right)
$$
and
$$
d(\left(\begin{array}{ccc}
\lambda&\sigma&0\\ \tau&\mu&0\\ 0&0&\nu
\end{array}\right))=
\left(\begin{array}{ccc}
0&0&a_{11}(\nu-\lambda)-a_{21}\sigma\\ 0&0&a_{21}(\nu-\mu)-a_{11}\tau\\ 0&0&0
\end{array}\right)
$$

We may study the homology of this differential graded algebra. Suppose that $R$ is an
integral domain.
If $a_{11}\neq 0\neq a_{21}$, then the degree $0$ cycles is
given by the set of matrices
$$
\left(\begin{array}{ccc}
\lambda&\frac{a_{11}}{a_{21}}(\nu-\lambda)&0\\ \frac{a_{21}}{a_{11}}(\nu-\mu)&\mu&0\\ 0&0&\nu
\end{array}\right)
$$
for $\lambda,\mu,\nu\in R$, such that all coefficients are in $R$.
This forms a $3$-dimensional affine variety. The boundaries in
degree $0$ forms a $2$-dimensional subvariety. Hence, the degree $0$ homology has an
$R$-torsion free part of rank $1$. Degree $-1$ cycles are trivial, since $R$ has no $0$ divisors.
The degree $1$ cycles is $R^2$, by definition. Clearly, the degree $1$ homology is
$\left(R/(a_{11}R+a_{21}R)\right)^2$.
\end{Example}

\section{Dg-algebras with respect to Semisimplicity and Noetherian ring theoretical properties}

\label{ringtheoryofdgalgebras}

As we shall need to consider semisimplicity of differential graded algebras, it makes sense to
consider differential graded versions of the Jacobson radicals. Further, we shall consider
semisimplicity of the homology algebra.

\subsection{Semisimplicity of the homology}

\begin{Example}
Consider Example~\ref{dgalgebraonthreebythree}. The special case $a_{11}=0\neq a_{21}$ produces an
interesting dg-algebra $A$, which is semisimple as an algebra.
Then
$$d(\left(\begin{array}{ccc}0&0 &0\\0 &0 &0\\x&y&0\end{array}\right)=a_{21}\cdot
\left(\begin{array}{ccc}0&0 &0\\x &y &0\\0&0&y\end{array}\right)$$
and
$$d(\left(\begin{array}{ccc}\lambda&\sigma &0\\\tau &\mu &0\\0&0&\nu\end{array}\right)=a_{21}\cdot
\left(\begin{array}{ccc}0&0 &-\sigma\\0 &0 &\nu-\mu\\0&0&0\end{array}\right).$$
Hence,
$$\ker(d)=\{\left(\begin{array}{ccc}\lambda&0 &s\\\tau &\mu &t\\0&0&\mu\end{array}\right)\;|\;u,v,s,t,x\in K\},$$
which is isomorphic to the semidirect product of the $2\times 2$ lower matrix algebra acting on its natural
$2$-dimensional representation: $\left(\begin{array}{c}K\\K\end{array}\right)\sdp\left(\begin{array}{cc}K&0\\K&K\end{array}\right).$
Clearly, this algebra is not semisimple.
However, $H(A)=K.$
%
\end{Example}

In this context is may be worth to mention the following criterion in this context.
We shall need the statement later.

\begin{Lemma}\label{semisimplekerd}
Let $K$ be a field and let $(A,d)$ be a finite dimensional differential graded $K$-algebra.
Then $(\ker(d),d)$ is a subalgebra of $(A,d)$. If $\ker(d)$ is semisimple, then
$H(A,d)$ is semisimple as well. Moreover, the restriction of the
epimorphism $\ker(d)\epi H(A,d)$ to $\ker(d)^\times$ yields a surjective group homomorphism
$\ker(d)^\times\lra H(A,d)^\times$.
\end{Lemma}

Proof. By Corollary~\ref{cyclesaresubalgebras} $\ker(d)$ is a graded subalgebra of $A$.
By hypothesis, $\ker(d)$ is a semisimple $K$-algebra.
Now, by Lemma~\ref{dgalgebrainducesalgebraonH}
there is an epimorphism $\ker(d)\stackrel{\pi}{\epi} H(A,d)$ of algebras,
and the preimage of a twosided
ideal in $H(A,d)$ is a twosided ideal of $\ker(d)$.
Therefore $H(A,d)\times\ker(\pi)=\ker(d)$ as algebras, and as a consequence any unit
$H(A,d)$ lifts to a unit of $\ker(d)$.
\dickebox

\begin{Rem} \label{complexmorphismsaresemisimple}
Let $(A,d)$ be a semisimple artinian dg-algebra and let $(L,d_L)$ and $(M,d_M)$ be
differential graded $(A,d)$-modules.
Then $\ker({(d_{Hom}})_0)=End_{(A,d)}(L,d)$ is not semisimple in general
since not every homomorphism $\alpha:(L,d_L)\ra (M,d_M)$
of complexes decomposes as $\alpha=\gamma\circ\delta$ for a split epimorphism
$\delta$ of complexes and a split monomorphism $\gamma$ of complexes.

We recall the easy argument.
We consider the special case when $A$ is concentrated in degree $0$ and $L$ is
finitely generated in each degree.

Since $A$ is semisimple, $L$ is a semisimple $A$-module.
We claim that $L=L_c\oplus H(L)$ is the direct sum of contractible $2$-term  complexes $L_c$
and its homology $H(L)$. Indeed,
since there is no non zero $A$-module homomorphism between two non isomorphic simple
$A$-modules, we may assume, multiplying by a specific central idempotent if necessary,
that $A$ actually is simple. By Morita equivalence we can therefore assume that
$A$ is a skew field.

{\red
Choose a basis of $\im(d_{-1})$, extend it to a basis of $\ker(d_0)$, and then the so-obtained basis
to a basis of $L_0$. For each basis element of $\im(d_{-1})$ choose a preimage in $L_{-1}$. The vector space
generated by these elements intersects $\im(d_{-2})$ in $0$, since $d^2=0$. Consider $\ker(d_{-1})$,
containing $\im(d_{-2})$. We may choose a basis of $\im(d_{-2})$, extend it to a basis
of $\ker(d_{-1})$, and complete the whole to a basis of $L_{-1}$
by the preimages chosen before. By downward induction we chose adapted bases in negative degrees.
For positive degrees, we proceed in the same way by upwards induction.
The bases elements chosen in the image of the differential, together with their
preimages give contractible $2$-term complexes. The rest combines to the homology.}


Consider the case of
$$L_1=\xymatrix{\dots 0\ar[r]&0\ar[r]&M\ar[r]^{id}&M\ar[r]&0\dots}$$
concentrated in degree $0$ and $1$, and
$$L_2=\xymatrix{\dots 0\ar[r]&N\ar[r]^{id}&N\ar[r]&0\ar[r]&0\dots}$$
concentrated in degree $-1$ and $0$. Moreover, suppose that $A=D$ is a skew field.
Then $Hom_{A}(L_1,L_2)=Hom_D(M,N)$ from the spaces in degree $0$, using that
the differentials do not impose any restriction. However,
$Hom_{A}(L_2,L_1)=0$ since the compatibility of the differentials impose
that the map in degree $0$ needs to be $0$.
%
%
%
\end{Rem}

{\red
\begin{Prop}\label{homologyissemisalwaysimple}
Let $(A,d)$ be a dg-algebra over a field $F$, such that $A$ is split simple artinian as an algebra.
Then $H(A,d)$ is simple as an algebra.
\end{Prop}

Proof.
By Proposition~\ref{progenerator} and Remark~\ref{strangehypothesisisunnecessary}
we see that $(A,d)\simeq (End_K^\bullet(C^\bullet),d_{Hom})$ for some finite dimensional
complex $C^\bullet$ of $F$-vector spaces. By Remark~\ref{complexmorphismsaresemisimple}
any bounded complex of vector spaces is isomorphic to $C^\bullet\simeq D^\bullet\oplus E^\bullet$
for $D^\bullet$ being contractible and $E^\bullet$ a $\Z$-graded vector space with zero differential.
Hence, we may replace $C^\bullet$ by this decomposition. In the homotopy category of
$F$-vector spaces we get $C^\bullet\simeq E^\bullet$.
Now, denoting by $K^b(F-mod)$ the homotopy category of bounded finite dimensional $F$-vector spaces,
\begin{eqnarray*}
H_n((A,d))&\simeq&H_n((End_K^\bullet(C^\bullet),d_{Hom}))\\
&\simeq&H_n((End_K^\bullet(D^\bullet\oplus E^\bullet),d_{Hom}))\\
&\simeq&Hom_{K^b(F-mod)}(D^\bullet\oplus E^\bullet,D^\bullet\oplus E^\bullet[n])\\
&\simeq&Hom_{K^b(F-mod)}( E^\bullet, E^\bullet[n])\\
&\simeq&Hom_{\textup{graded $F$-modules}}( E^\bullet, E^\bullet[n])
\end{eqnarray*}
Therefore,  $H(A,d)$ is a matrix algebra over $F$,
and the grading on the homology gives a grading on the matrix algebra. \dickebox

\begin{Rem}
The hypothesis that $A$ is split is not necessary since actually $F$ may well be a skew field
in the proof of Proposition~\ref{progenerator}.
\end{Rem}
}

\begin{Rem}\label{semisimpleendorings}
Consider the situation of Remark~\ref{complexmorphismsaresemisimple}.
Then we have a ring homomorphism
$$\xymatrix{End_A^\bullet(L)\ar[r]&H_0(End_A^\bullet(L))=End_{gr-A}(H(L))}$$
from the endomorphism complex of $L$ to the ring of graded $A$-linear
homomorphisms of $H(L)$.
Since $A$ is semisimple, and since $H(L)$ is a semisimple $A$-module,
$End_{gr-A}(H(L))$ is a semisimple algebra.
Moreover, using Lemma~\ref{Hombullet}, we get that $(End_A^\bullet(L),d_{Hom})$ is a differential
graded algebra, mapping to a semisimple algebra $End_A(H(L))$.
\end{Rem}

\subsection{Differential graded radicals of differential graded algebras}

%

We shall study the ring theory of differential graded algebras.

\begin{Lemma}\label{differentialgradedquotient}
Let $(\Lambda,d)$ be a differential graded algebra, let $(M,\delta)$ be a differential graded
$(\Lambda,d)$-module with a differential graded submodule $(N,\delta)$. Then there is a
differential $\overline\delta$ on $M/N$ given by $\ol\delta(m+N):=\delta(m)$ for all $m\in M$.
\end{Lemma}

Proof. $M/N$ is a $\Lambda$-module by general ring theory.
Since $\delta(N)\subseteq N$, the map $\ol\delta$ is well-defined. $\ol\delta$ is
clearly additive since $\delta$ is additive. Since $\delta^2=0$, also $\ol\delta^2=0$.
Since $\delta$ satisfies Leibniz formula, so is $\ol\delta$. \dickebox

\begin{Lemma}\label{dgmaximalideals}
Let $(\Lambda,d)$ be a differential graded algebra and let $(I,d)$ be a differential
graded ideal of $(\Lambda,d)$.
Then there is a differential graded ideal $(M,d)$ in $(\Lambda,d)$, {\red which is maximal in the
partial ordered set of differential graded ideals different from $(\Lambda,d)$, and} containing $(I,d)$.
\end{Lemma}

Proof. An analogue of the usual proof in the classical situation works. We use Zorn's lemma.
Let $\mathcal X$ be the set of differential graded ideals
of $(\Lambda,d)$ containing $(I,d)$ but
different from $(\Lambda,d)$.
Then $(I,d)\in  \mathcal X$, and hence $\mathcal X$
is not empty. Let ${\mathcal L}$ be a totally ordered subset
of $\mathcal X$ and let $J:=\bigcup_{L\in{\mathcal L}}L$. Then $J$ is an ideal by the
classical argument. Moreover, it is differential graded since
if $x\in J$, there is $L_x\in\mathcal L$
with $x\in L_x$. Then $d(x)\in L_x\subseteq J$. Further the
Leibniz formula holds in $\Lambda$,
whence also in $J$. If $1\in J$, then there is $L_1\in\mathcal L$ such that $1\in L_1$,
and hence $L_1$ is
not in $\mathcal X$, which provides a contradiction since
${\mathcal L}\subseteq{\mathcal X}$. \dickebox

\begin{Lemma}\label{dgfield}
Let $(\Lambda,d)$ be a differential graded algebra and let $(I,d)$ be a differential
graded ideal of $(\Lambda,d)$. Then $(I,d)$ is a differential graded ideal, {\red  maximal in the partial ordered
set of proper differential graded ideals},
if and only if $(\Lambda/I,\ol d)$
does not contain non trivial differential graded ideals, for $\ol d$ the differential
induced on $\Lambda/I$
by $d$.
\end{Lemma}

Proof. Let $x\in\Lambda$ and $y\in I$. Then by Lemma~\ref{differentialgradedquotient} for
$\ol d(x+I):=d(x)+I$ defines a differential on $\Lambda/I$.
Consider the zero ideal in $(\Lambda/I,\ol d)$. By Lemma~\ref{dgmaximalideals}
this is contained in a maximal {\red (in the above sense)}
differential graded ideal $(\ol M,\ol d)$. Let $M$ be the preimage of
$\ol M$ in $\Lambda$. We claim that $(M,d)$ is a dg-ideal.
Since $\ol M$ is an ideal, by classical ring theory $M$ is an ideal as well.
For $m\in M$ let $d(m)+I=\ol d(m+I)\in \ol M$ since $(\ol M,\ol d)$ is a
dg-ideal. Hence $d(m)\in M$ and by consequence $(M,d)$ is a differential graded ideal.
If there is no non trivial differential graded ideal in $\Lambda/I$, then $M=I$ and $I$ is a
maximal {\red (in the above sense)}
differential graded ideal. If $\Lambda/I$ admits non trivial differential graded
ideals, then $I\neq M$, and $I$ is not maximal.
This shows the statement. \dickebox

\begin{Rem}
Note that the statement of Lemma~\ref{dgfield} is completely formal and
holds for left- or right- or twosided ideals.
\end{Rem}

{\red
\begin{Lemma}\label{mfinitelygenerateddg}
Let $(A,d)$ be a differential graded ring and let $(M,\delta)$ be a non zero differential graded
$(A,d)$-module. If $(M,\delta)$ is finitely generated as differential graded module, then
there is a differential graded submodule $(N,\delta)$ of $(M,\delta)$, maximal with respect to
the lattice of proper differential graded submodules, and such that
$(M/N,\ol\delta)$ is a dg-simple dg-module.
\end{Lemma}

Proof. Let $(M,\delta)$ be finitely generated
and denote by $\{g_1,\dots,g_n\}$ a set of generators.
We need to show that $M$ contains a proper dg-submodule $(N,\delta)$ which is maximal as dg-submodule.
Let $\mathcal X$ be the set of dg-submodules $(L,\delta)$ of $(M,\delta)$ such that $L\neq M$.
This set is not empty since $0$ is trivially a dg-submodule.
Let $\mathcal L$ be a non empty totally ordered subset of $\mathcal X$.
Then, analogous to the proof of Lemma~\ref{dgmaximalideals}
$$V:=\bigcup_{L\in{\mathcal L}}L$$
is a dg-submodule of $(M,\delta)$. We need to show that $V\neq M$. Else, $\{g_1,\dots,g_n\}\subseteq M= V$
and hence for any $i\in\{1,\dots,n\}$ there is $L_i\in{\mathcal L}$ with $g_i\in L_i$.
Since $\mathcal L$ is totally ordered, there is $m\in\{1,\dots,n\}$ such that $L_i\subseteq L_m$
for all $i\in\{1,\dots,n\}$. But then $M=L_m$
since $L_m$ contains all generators $\{g_1,\dots,g_n\}$. This is a contradiction to $L_m\in{\mathcal X}$.
By Zorn's lemma there is a maximal element in $\mathcal X$, say $(N,\delta)$.
We claim that $N$ is maximal with respect to the lattice of dg-submodules. Else,
let $(D,\delta)$ be a strictly bigger proper dg-submodule. Then $(D,\delta)$ is
superior to $(N,\delta)$ in $\mathcal X$, a contradiction.
Since $(N,\delta)$ is maximal as a dg-submodule, $(M/N,\ol\delta)$ is dg-simple, since else
consider a proper quotient of $(M/N,\ol\delta)$, and the kernel of the composite quotient map
from $M$ to the quotient would be strictly larger than $(N,\delta)$.
We proved the statement. \dickebox
}

%
%
%

\begin{Lemma}\label{prinsipaldgideals}
Let $(\Lambda,d)$ be a differential graded algebra {\red and let $(M,\delta)$ be a differential
graded left (resp. right) module.} Then for any {\red homogeneous} $a\in ker(\delta)$
we have that $\Lambda\cdot a$ (resp. $a\cdot\Lambda$) is a differential graded $(\Lambda,d)$ left
(resp. right)  submodule.
\end{Lemma}

Proof. Indeed, $\delta(\lambda\cdot a)=d(\lambda)\cdot a\pm \lambda \delta(a)=d(\lambda)\cdot a$
and likewise $\delta(a\cdot\lambda )=\delta(a)\lambda\pm a\cdot d(\lambda)=a\cdot d(\lambda)$.
\dickebox

\begin{Rem}\label{principaldgideal}
If the {\red homogeneous} $a$ is not (necessarily) a cycle, i.e. $a\not\in ker(d)$, then
$\Lambda \cdot a+\Lambda\cdot d(a)$ is a differential graded ideal.
Similarly, if $(M,d_M)$ is a dg-bimodule, and if $m\in\ker(d_M)$, then
$\Lambda\cdot m+\Lambda\cdot d_M(m)$ is a left dg-submodule of $(M,d_M)$.
\end{Rem}

\begin{Cor}\label{maximaldgidealandinvertiblesinquotient}
Let $(\Lambda,d)$ be a differential graded algebra and let $(I,d)$ be a
left (resp. right) differential graded ideal, {\red maximal in the lattice of proper left (resp. right) dg-ideals}.
If $a\in (ker(d)+I)$ {\red is homogeneous},  and if $a\not\in I$, then $a+I$ is a left (resp.
right) generator of the module $\Lambda/I$.
\end{Cor}

Proof. If $a\in ker(d)+I$, but $a\not\in I$,
then $a+I\in ker(\ol d)$ and then there is $x\in I$ such that $a+x\in\ker(d)$ and hence
$\Lambda\cdot (a+x)$
is a differential graded
ideal of $\Lambda$. Hence, also $I+\Lambda(a+x)=I+\Lambda a$ is a
differential graded ideal of $\Lambda$.
If $a\not\in I$, then by the maximality of $(I,d)$ we have that $\Lambda=I+\Lambda a$.
Therefore $a+I$ has a left inverse in $\Lambda/I$. The statement for the right ideals is analogous. \dickebox

%
%

\begin{Def} \label{semisimpledef}
Let $(\Lambda,d)$ be a differential graded algebra.
A non zero differential graded $(\Lambda,d)$-module $(M,\delta)$
is called {\em differential graded simple} if there is no differential graded submodule
of $(M,\delta)$ different from $(M,\delta)$ and $0$.
A differential graded $(\Lambda,d)$-module $(M,\delta)$ is called {\em differential
graded semisimple} if
$(M,\delta)$ is the direct sum of simple differential graded $(\Lambda,d)$-modules.
A {\em differential graded algebra $(\Lambda,d)$ is called  semisimple} if
$(\Lambda,d)$ as left $(\Lambda,d)$-module is semisimple.
\end{Def}

{\red
\begin{Rem}
We should remind the reader of different concepts of semisimplicity.
In Definition~\ref{semisimpledef} we consider a dg-module as being semisimple if it
is a direct sum of simple dg-modules, and a dg-module is dg-simple if it does not contain a
proper dg-submodule.

In contrast to this is the concept of a dg-module being semisimple if all dg-submodules
admit a dg-complement. This concept is the one studied by Aldrich and
Garcia-Rozas~\cite{Tempest-Garcia-Rochas}, and they obtain the complete classification
as displayed in Theorem~\ref{GarciaRozastheorem}.

For differential graded modules over differential graded algebras these two concepts
differ.

Recall that Lam~\cite{Lamfirstcourse} give two notions for these properties.
He says (cf \cite[I (2.1) Definition]{Lamfirstcourse}) that a module over an algebra
is semisimple if all submodules admit a complement. He calls an algebra
Jacobson-semisimple (or J-semisimple) if its Jacobson
radical is $0$ (cf \cite[II (4.7) Definition]{Lamfirstcourse}).
\end{Rem}
}

\begin{Example} \label{dualnumbersexample}
Let $K$ be a field
and let $A=K[X]/X^2$ for $deg(X)=-1$. Then $d(1)=0$ and $d(X)=1$
give a structure of differential graded algebra to $A$. Indeed,
$(a+bX)(c+dX)=ab+(bc+ad)X$ and $$d((a+bX)(c+dX))=bc+ad$$
and
$$
d(a+bX)\cdot (c+dX)+(a-bX)\cdot d(c+dX)=b(c+dX)+(a-bX)d=bc+ad.
$$
Further, $0=d(0)=d(X^2)=d(X)\cdot X-X\cdot d(X)=X-X=0$.
The algebra $A$ is simple as differential graded algebra {\red in the sense that
there is no non trivial twosided differential graded
ideal (cf Definition~\ref{dgsimpledgalgebradef} below)}. Indeed, there is only one non
trivial algebra ideal, namely $XK[X]$. However, this is not a differential graded
ideal, since $d(X)=1$. Further, the only simple differential graded modules are the regular one
{\red and its shifts in degree},
which are not concentrated in a single degree. Note that this is the smallest example
of algebras mentioned in the second item of Theorem~\ref{GarciaRozastheorem}
studied by Aldrich and Garcia Rozas \cite{Tempest-Garcia-Rochas}.
Of course, $A$ is not a simple algebra.
\end{Example}

\begin{Cor}\label{quotientsmodmaxaresimple}
Let $(\Lambda,d)$ be a differential graded algebra and let $(I,d)$ be differential graded
left ideal of $(\Lambda,d)$ {\red which is maximal in the lattice of proper dg-ideals}.
Then $(\Lambda/I,\ol d)$ is a {\red dg-}simple differential graded $(\Lambda,d)$-module.
\end{Cor}

Proof. This is a direct consequence of Lemma~\ref{dgfield}. \dickebox

\begin{Lemma}\label{simplesarequotients}
Let $(\Lambda,d)$ be a differential graded algebra and let $(S,\delta)$ be a {\red dg-}simple
differential graded $(\Lambda,d)$-module. Then for any non zero homogeneous $u\in\ker(\delta)$ there is a
surjective homomorphism $(\Lambda,d)[-|u|]\stackrel{\pi}\lra (S,\delta)$ given by $\pi(\lambda)=\lambda u$.
Further, there is always a non zero homogeneous $u\in\ker(\delta)$, and $\ker(\pi)$ is a
differential graded ideal {\red which is maximal in the set of proper dg-ideals}.
\end{Lemma}

Proof. Due to the shift, $\pi$ is a homomorphism of graded modules.
We need to show that $\pi$ is compatible with the differential. Let $\lambda\in\Lambda$ be homogeneous.
Then
$$\pi(d(\lambda))=d(\lambda)u=d(\lambda)u+(-1)^{|\lambda|}\lambda\delta(u)=\delta(\lambda u)=\delta(\pi(\lambda)),$$
where the third equation is Leibniz formula, and where the second equation holds
since $u\in\ker\delta$.
Using Lemma~\ref{prinsipaldgideals} we see that
$\im(\pi)$ is a non trivial dg submodule of $S$, and by simplicity of $S$, we see that
$\pi$ is surjective. Note that $\delta^2=0$ shows that $\ker(\delta)$ is not zero.
Using Lemma~\ref{dgfield} we see that $\ker(\pi)$ is a maximal dg-ideal.
\dickebox

\begin{Lemma}\label{semisimpledgimpliescomplemented}
Let $(\Lambda,d)$ be a differential graded algebra, let $(M,\delta)$ be a differential graded
$(\Lambda,d)$-module. Suppose that  there are {\red dg-}simple differential graded $(\Lambda,d)$-left modules
$(S_1,\delta_1)$, ..., $(S_n,\delta_n)$ such that
$$(M,\delta)\simeq (S_1,\delta_1)\oplus\dots\oplus (S_n,\delta_n)$$
as differential graded left modules.
Then for any differential graded submodule $(N,\delta)$ of $(M,\delta)$ there is a differential graded
submodule $(L,\delta)$ of $(M,\delta)$ such that $M=N\oplus L$.
\end{Lemma}

Proof.
Consider the family of sets $\{i_1,\dots,i_k\}\subseteq\{1,\dots,n\}$ such that
$\left(\bigoplus_{j=1}^kS_{i_j}\right)\cap N=\{0\}$. Then taking a maximal one $P=\{i_1,\dots,i_k\}$
 of these subsets, and
put $L:=\bigoplus_{j=1}^kS_{i_j}$. Then $L+N=L\oplus N$ is a differential graded
$(\Lambda,d)$-module. Hence we
get that $(S_\ell+L)\cap N\neq\{0\}$ for any $\ell\not\in P$.
Therefore $(L\oplus N)\cap S_\ell\neq\{0\}$ for any $\ell\not\in P$. Since $S_\ell$
is a differential graded simple $(\Lambda,d)$-module, and since
$L\oplus N$ is a differential graded module, $S_\ell\subseteq L\oplus N$ for
any $\ell\not\in P$. But this shows that $L\oplus N\supseteq\bigoplus_{i=1}^nS_i=M$.
We proved the statement. \dickebox

\begin{Def}\label{dgsimpledgalgebradef}
Let $K$ be a field and let $(A,d)$ be a differential graded $K$-algebra.
The dg-algebra $A$ is called {\em differential graded simple} if the only
differential graded twosided $(A,d)$-ideals are $0$ and $(A,d)$.
\end{Def}

\begin{Rem}
Recall from the second item of Theorem~\ref{GarciaRozastheorem} that every differential graded $(A,d)$-module
is semisimple if and only if $(A,d)$ is {\red acyclic} and $\ker(d)$ is a semisimple algebra.
Moreover, then $A=\ker(d)\oplus z\cdot\ker(d)$ for $z\in d^{-1}(1)$.

Consider for example the case of an algebra concentrated in degree $0$, such as the field $K$. Then
there is no non trivial twosided ideal. However, this dg-algebra does not satisfy the condition
of Theorem~\ref{GarciaRozastheorem}.
The classical result that simple artinian algebras are semisimple cannot be transposed
to the dg-situation.
\end{Rem}

\begin{Prop} \label{simpledgalgebras}
Let $K$ be a field and let $(A,d)$ be a {\red dg-simple}
differential graded algebra.
Suppose that the regular differential graded left $(A,d)$ module and the
regular differential graded right $(A,d)$ module are {\red dg-}artinian.

Denote by $\textup{cone}(id_K)$ the cone of the identity map of the stalk complex $K$ in $D^b(K-mod)$.

Then, as differential graded left $(A,d)$-module, the regular module
$(A,d)$ is the direct sum of shifts of
copies of one specific simple differential graded
modules $(S,\delta)$ and shifted copies of modules of the form total complex of the tensor product
of $(S,\delta)$ by $\textup{cone}(id_K)$, i.e.
$tot\left((S,\delta)\otimes_K(0\lra K\stackrel{id}\lra K\lra 0)\right)$.
%
\end{Prop}

Proof.
As $A$ is dg-artinian there is a minimal non zero differential graded left ideal $(S,d)$
and a minimal non zero differential graded right ideal $(T,d)$.
Then consider the $(A,d)-(A,d)$-bimodule $X:=(S,d)\otimes_K(T,d)$.
The multiplication map
$$A\otimes_KA\lra A$$
is a chain map (replacing $A\otimes_KA$ by $tot(A\otimes_KA)$) by the Leibniz rule and
maps $X$ to a bimodule, whence a twosided ideal of $A$.
Hence the image is a differential graded twosided ideal of $(A,d)$. Since there is no non trivial
twosided ideal, this image is $A$. Now, consider $X$ as differential graded left module.
Since only the $K$-structure of $(T,d)$ is of importance then, and since complexes
of $K$-modules are a direct sum of shifted copies of stalk complexes
$K$ and shifted copies of complexes $\textup{cone}(id_K)$.
This shows the statement. \dickebox

\begin{Rem}\label{exampleofdgsimpleoverdgsimple}
We should note that the short exact sequence
$$
\xymatrix{
&0\ar[r]&K\ar@{^{(}->}[d]\ar[r]&0\\
0\ar[r]&K\ar[r]\ar@{->>}[d]&K\ar[r]&0\\
0\ar[r]&K\ar[r]&0
}
$$
of differential graded $K$-modules yields a
short exact sequence
$$0\lra (S,\delta)\lra tot\left((S,\delta)\otimes_K(0\lra K\stackrel{id}\lra K\lra 0)\right)\lra
(S,\delta)[1]\lra 0$$
of $(A,d)$-left modules.
{\red Hence, this tensor product is not a dg-simple dg-module.}
\end{Rem}

\bigskip

For any differential graded left $(\Lambda,d)$-module $(M,\delta)$ denote the {\em
annihilator of $(M,\delta)$} by
$$\ann(M,\delta):=\{\lambda\in\Lambda\;|\;\lambda m=0\;\;\forall m\in M\}. $$

\begin{Lemma}\label{annistwosideddg}
We get that $\ann(M,\delta)$ is a differential graded twosided ideal of $(\Lambda,d)$.
\end{Lemma}

Proof.
Indeed, $\ann(M,\delta)$ is a twosided {\red graded} ideal by general ring theory.
Further, if $\lambda\in \ann(M,\delta)$ is homogeneous, then $\lambda m=0$ for any $m\in M$, and hence
$0=\delta(\lambda m)=d(\lambda)m\pm \lambda \delta(m)=d(\lambda)m$ since $\delta(m)\in M$
and $\lambda\in \ann(M,\delta)$. This show that $d(\lambda)m=0$, and hence $d(\lambda)\in \ann(M,\delta)$.
\dickebox

\begin{Def}
Let $(\Lambda,d)$ be a differential graded algebra.
\begin{itemize}
\item Then the {\em differential graded left radical}
$\dgrad_\ell(\Lambda,d)$ is the intersection of all dg-maximal differential graded
left ideals of $(\Lambda,d)$.
\item Likewise the {\em differential graded right radical} $\dgrad_r(\Lambda,d)$
is the intersection of all dg-maximal differential graded right ideals of $(\Lambda,d)$.
\item The {\em differential graded twosided radical}
$\dgrad_2(\Lambda,d)$ is the intersection of all of
the annihilators of {\red dg-}simple differential graded left modules.
\end{itemize}
\end{Def}

\begin{Lemma} (differential graded Nakayama's Lemma)\label{dgNakayama}
Let $(\Lambda,d)$ be a differential graded algebra.
{\red
Let $M$ and $N$ be  differential graded $\Lambda$-modules
{\red such that $N\leq M$}. Assume that $(M,\delta)$ is finitely generated as a dg-module.
Then $\dgrad_2 (\Lambda)M=M$ implies
$M=0$, and $N+\dgrad_2 (\Lambda)M=M$ implies $M=N$.}
\end{Lemma}

Proof. {\red
We  suppose that $(M,\delta)$ is finitely generated as a dg-module.
Again, in case $M=0$ the statement is trivial.
By Lemma~\ref{mfinitelygenerateddg} the dg-module $(M,\delta)$ contains a proper
maximal dg-submodule $(D,\delta)$ and $(M/D,\ol\delta)$ is a non trivial dg-simple dg-module.
By definition of $\dgrad_2(\Lambda,d)$ we get $\dgrad_2(\Lambda,d)\cdot M/D=0$, whence
$\dgrad_2(\Lambda,d)\cdot M\subseteq D$. This proves that only $M=0$ satisfies
$\dgrad_2(\Lambda,d)\cdot M=M$.} As for the second part,
$N+\dgrad_2 (\Lambda)M=M$ implies $\dgrad_2 (\Lambda)(M/N)=M/N$. The first statement then
implies $M/N=0$, and therefore $M=N$.
Hence the result follows.
\dickebox

{\red \begin{Rem}
The statement of the dg-Nakayama Lemma for finitely generated dg-modules $M$
was suggested by the referee. I am very grateful for this suggestion.
\end{Rem}}

{\red
\begin{Lemma}\label{dgnoetherianfg}
A dg-module is dg-Noetherian
if and only if all its dg-submodules are finitely generated as dg-modules.
\end{Lemma}

Proof.
The classical proof of this fact transposes to the dg-case, after slight modifications.
Let us verify the details.

If $(M,\delta)$ is dg-Noetherian, and the dg-submodule $(N,\delta)$ is not finitely
generated, then take a non zero $g_1\in\ker(\delta)\cap N$. Such an element exists since
if $g_1\in N\setminus\ker(\delta)$, then $\delta(g_1)\in\ker(\delta)\cap N$.
Using Lemma~\ref{prinsipaldgideals}
we get that $A\cdot g_1$ is a dg-submodule of $N$. Suppose we constructed $g_1,\cdots,g_{n-1}$ such that $\sum_{i=1}^jAg_i$
is a dg-submodule for all $j$ and $\sum_{i=1}^jAg_i\lneq\sum_{i=1}^{j+1}Ag_i$ for all $j$.
Since $(N,\delta)$ is not finitely generated, there is $g_{n}\in N\setminus \sum_{i=1}^{n-1}Ag_i$.
Consider the dg-module $(N_1,\ol\delta):=(N/\sum_{i=1}^{n-1}Ag_i,\ol\delta)$.
If $\ol\delta(g_n+\sum_{i=1}^{n-1}Ag_i)\neq 0$, then replace $g_n$ by $\delta(g_n)$, which still is not in $\sum_{i=1}^{n-1}Ag_i$, and
the image in $N/\sum_{i=1}^{n-1}Ag_i$ of this replaced element is in $\ker(\ol\delta)$.
If $\ol\delta(g_n+\sum_{i=1}^{n-1}Ag_i)=0$ there is $h_n\in \sum_{i=1}^{n-1}Ag_i$ with $\delta(g_n)=h_n$. But
then $\sum_{i=1}^{n}Ag_i$ is a dg-submodule of $N$. By induction
we obtain a strictly increasing sequence of dg-submodules
$$Ag_1\subsetneq Ag_1+Ag_2\subsetneq  Ag_1+Ag_2+Ag_3\subsetneq\cdots$$
of $(N,\delta)\subseteq (M,\delta)$. This contradicts the hypothesis that $(M,\delta)$ is dg-Noetherian.

If all dg-submodules of $(M,\delta)$ are finitely generated, then let
$$N_1\subseteq N_2\subseteq N_3\subseteq\cdots\subseteq M$$
be a deceasing sequence of dg-submodules. Put $L:=\bigcup_{i=1}^\infty N_i$. This is a
dg-submodule of $M$, and hence is finitely generated, by $\{g_1,\dots,g_n\}$ say, by hypothesis.
For each $i\in\{1,\dots,n\}$ there is $j(i)$ such that $g_i\in N_{j(i)}$. But then, for $m=\max(j(1),\dots,j(n))$
we have $\{g_1,\dots,g_n\}\subseteq N_m$, and therefore $N_m=N_{m+k}$ for all $k\geq 0$.
Therefore $M$ is dg-Noetherian.
\dickebox

\begin{Cor}
Let $M$ and $N$ be differential graded $\Lambda$-modules
{\red such that $N\leq M$}. If $M$ is {\red dg-}Noetherian,
then $\dgrad_2 (\Lambda)M=M$ implies
$M=0$, and $N+\dgrad_2 (\Lambda)M=M$ implies $M=N$.
\end{Cor}

Proof. By Lemma~\ref{dgnoetherianfg} the hypothesis of Lemma~\ref{dgNakayama} are satisfied. \dickebox

\begin{Rem}
Note that we need Zorn's lemma in the proof of Lemma~\ref{dgNakayama}. If $M$ is
dg-Noetherian and dg-artinian, then a proof without Zorn's lemma is possible by a short
induction on the composition length.

If $M=0$, the (first) statement is trivial. Else, the hypothesis on $M$ implies that
there is a dg-simple dg-module, namely a dg-simple dg-submodule of $M$.
If $\dgrad_2 (\Lambda)M=M$, then comparing the composition lengths, one needs to have $M=0$.
\end{Rem}
}

{\red \begin{Rem}
There is a parallel work on Nakayama's Lemma for differential graded algebras by
Goodbody~\cite{Goodbody}, based on the work of Orlov~\cite{Orlov1}. These results are
very different from ours, and in particular Goodbody uses that the algebra is finite
dimensional over some field.
I discovered
both articles after having submitted this manuscript.
\end{Rem}}

\begin{Prop} \label{dgradistwosidedideal}
Let $(\Lambda,d)$ be a differential graded algebra, and suppose that $(\Lambda,d)$ is artinian
as differential graded left and as a differential graded right module.
Then $$\dgrad_\ell(\Lambda,d)\cap \dgrad_r(\Lambda,d)\supseteq \dgrad_2(\Lambda,d).$$
\end{Prop}

Proof. We shall prove that $\dgrad_2(\Lambda,d)\subseteq\dgrad_\ell(\Lambda,d)$. The analogous statement
for $\dgrad_r(\Lambda,d)$ is symmetric.
By {\red Corollary~\ref{quotientsmodmaxaresimple}} every simple differential graded $(\Lambda,d)$-module is
isomorphic to one of the form
$(\Lambda/M,\ol d)$ for a maximal differential graded ideal $(M,d)$ of $(\Lambda,d)$.
{\red By definition of $\dgrad_2(\Lambda,d)$ the dg-simple dg-module $(\Lambda/M,\ol d)$
is annihilated by $\dgrad_2(\Lambda,d)$. Hence $\dgrad_2(\Lambda,d)\subseteq M$.
This shows $\dgrad_2(\Lambda,d)\subseteq\dgrad_\ell(\Lambda,d)$.}
This shows the proposition. \dickebox

{\red \begin{Rem}
The hypothesis that $(\Lambda,d)$ is dg-artinian is only used to show that
there are dg-simple dg-left modules. Just like in the classical situation, an easy use of
Zorn's lemma shows that it is sufficient to ask that $(\Lambda,d)$ is finitely generated
as dg-left module.
\end{Rem}}

{\red
\begin{Example}
Let $K$ be a field.
Recall from Example~\ref{Examplediffgradedalgebra}
that the graded endomorphism algebra $\left(\begin{array}{cc}K&K\\K&K\end{array}\right)$
of the complex $\textup{cone}(\textup{id}_K)$  is a dg-algebra.
It is not difficult to verify that $\left(\begin{array}{cc}0&K\\0&K\end{array}\right)=:I$
is the only non trivial differential graded left ideal of $A$. Hence $\dgrad_\ell(A)=I$
and likewise the first row is the only non trivial differential graded right ideal.
Hence $\dgrad_\ell(A)\cap\dgrad_r(A)=\left(\begin{array}{cc}0&K\\0&0\end{array}\right)$
is not an ideal at all. Further it strictly contains $\dgrad_2(A)=0$.
\end{Example}
}

\begin{Prop} \label{moddgrad}
Let $(\Lambda,d)$ be a differential graded algebra. Suppose that $(\Lambda,d)$ is
dg-Noetherian and dg-artinian, i.e.
Noetherian and artinian as differential graded left module over itself.
Then $\Lambda/\dgrad_\ell (\Lambda)$ is {\red a direct sum of dg-simple dg-modules
over $(\Lambda,d)$}. In particular,
$\dgrad_\ell(\Lambda)=0$ if and only if any {\red finitely generated }
differential graded $(\Lambda,d)$-module is {\red isomorphic}
to a direct sum of {\red dg-}simple differential graded
$(\Lambda,d)$-modules.
\end{Prop}

Proof. 
{\red In general there are infinitely many dg-maximal dg-left ideals.}
We shall show that there are only finitely many dg-maximal differential graded left
ideals $I_1,\dots,I_n$  of $\Lambda$ such that $I_j\not\supseteq \bigcap_{i=1}^{j-1}I_i$
{\red for all $j$.}
Take a maximal differential graded left ideal $(I_1,d)$. If there is no other
maximal differential graded left ideal $(I_2,d)$, we are done.
Else we consider $(I_1\cap I_2,d)$. If there is no maximal differential graded left ideal $I_3$
with $(I_3,d)\not\supseteq I_1\cap I_2$, we are done. Else we consider $(I_1\cap I_2\cap I_3,d)$.
Inductively we obtain a chain of strictly decreasing differential graded left ideals.
Since $\Lambda$ is assumed to be dg-artinian, this chain is finite. Hence for any set of maximal
differential graded ideals there is a finite subset
of dg-maximal differential graded ideals $I_1,\dots,I_n$ such that $\bigcap_{i=1}^nI_i$
is contained in any dg-maximal differential graded left ideal of $(\Lambda,d)$.
Hence
$$I_1\cap\dots\cap I_n=\dgrad_\ell(\Lambda)=\bigcap_{I\textup{ {\red dg-}maximal dg-ideal}}I.$$

Then there is a homomorphism
$$\Lambda\stackrel{\rho}\lra \Lambda/I_1\times\dots\times\Lambda/I_n$$
with kernel $\dgrad_\ell (\Lambda)$.
Hence $\rho$ induces a monomorphism
$$\Lambda/\dgrad_\ell (\Lambda)\stackrel{\ol\rho}{\hookrightarrow}  \Lambda/I_1\times\dots\times\Lambda/I_n.$$
Lemma~\ref{semisimpledgimpliescomplemented} then shows that there is a differential graded
module $M$ such that
$$\Lambda/\dgrad_\ell (\Lambda)\oplus M\simeq \Lambda/I_1\times\dots\times\Lambda/I_n.$$
Since $\Lambda$ is dg-Noetherian and dg-artinian, the Krull-Schmidt theorem
applies and hence, renumbering if necessary,
there is $k$ such that
$$\Lambda/\dgrad_\ell (\Lambda)\simeq \Lambda/I_1\times\dots\times\Lambda/I_k\textup{
and
}M\simeq \Lambda/I_{k+1}\times\dots\times\Lambda/I_n.$$
However, by construction of the sequence of ideals $I_1,\dots,I_n$  we have $k=n$.
Hence $\ol\rho$ is an isomorphism.

If $\Lambda$ is a direct sum of {\red dg-}simple differential graded left modules,
$\Lambda=S_1\times\dots\times S_n$, and hence
a maximal dg-ideal is $S_1\times\dots\times S_{i-1}\times S_{i+1}\times\dots S_n$. Therefore the
intersection of all dg-maximal dg-ideals is $0$.

{\red Let $(X,\delta)$ be a finitely generated differential graded $(\Lambda,d)$-module.
Then
there is an epimorphism $\bigoplus_{j=1}^n\Lambda[n_j]\lra X$. By
Lemma~\ref{semisimpledgimpliescomplemented}, the kernel is a direct factor of $\bigoplus_{j=1}^n\Lambda[n_j]$
and this shows the statement.}

This shows the proposition. \dickebox

\medskip

{\red
We shall now characterise $\dgrad_2(A,d)$ by universal properties. For this purpose we need
to recall some classical concepts.

Recall that an algebra $A$ is called left primitive if $A$ allows a faithful simple left $A$-module.
Simple rings allowing a simple module are trivially primitive, and the converse is true for
artinian algebras. In general primitive algebras are not simple. The standard
example is the endomorphism algebra (as a vector space) of a vector space of countable dimension.
The set of endomorphisms with finite dimensional image is a twosided ideal. However, the endomorphism
algebra is primitive. We further mention that in general a simple (hence primitive) algebra
can have many non isomorphic simple modules. An example is the Weyl algebra over an
algebraically closed field of characteristic $0$
(see Block~\cite{Block,Block2}). There are algebras which are left primitive, but not right primitive.
Since these properties hold for ordinary algebras, and since an algebra is a
dg-algebra with trivial grading and differential $0$, we cannot expect better properties
than these in the dg-version.

Recall further that for algebras $A_i$ indexed by elements of a set $i\in I$, a subdirect product is
defined to be a subalgebra $S$ of $\prod_{i\in I}A_i$ such that the composition of the injection
of $S$ into the product, followed by the projection onto $A_i$ is onto for any $i\in I$.
In particular, the intersection of the kernels of these composition of injection and projections is $0$
in $S$.

We come to the straight forward corresponding dg-concept.

\begin{Def}
A dg-algebra $(A,d)$ is called {\em dg-primitive} if there is a faithful dg-simple
differential graded left module $(S,\delta)$ over $(A,d)$. A twosided dg-ideal $(I,d)$
of $(A,d)$
is called {\em left dg-primitive} if $(A/I,\ol d)$ is a left dg-primitive algebra
\end{Def}

\begin{Lemma}\label{lambdamodannhasonlyonesimple}
Let $(\Lambda,d)$ be a
differential graded algebra and let $(S,\delta)$ be a dg-simple
differential graded module. Then, $\ann(S,\delta)$ is a  dg-primitive
differential graded twosided ideal.
\end{Lemma}

Proof. The fact that $\ann(S,\delta)$ is a twosided dg-ideal was already established.
Since $(S,\delta)$ is a simple $(\Lambda,d)$-module, it is still a
simple $(\Lambda/\ann(S,\delta))$-module. By definition, it is a faithful
$(\Lambda/\ann(S,\delta))$-module. \dickebox}

{\red
\begin{Rem}
Note that we always have $\ann(M,\delta)=\ann((M,\delta)[1])$, and if $(M,\delta_M)\simeq(N,\delta_N)$, then
$\ann(M,\delta_M)=\ann(N,\delta_N)$.
\end{Rem}

\begin{Prop}  \label{moddgradisprodsimples}
Let $(\Lambda,d)$ be a differential graded algebra and suppose that $(\Lambda,d)$ allows
a dg-simple dg-module. {\red Then
$\dgrad_2(\Lambda,d)$ is a twosided
differential graded ideal.}  Then $\dgrad_2(\Lambda,d)$ is
the smallest twosided differential graded ideal $I$ such that
$(\Lambda/I,\ol d)$ is a subdirect product of dg-primitive differential graded algebras.
\end{Prop}

Proof.
The fact that $\dgrad_2(\Lambda,d)$ is a twosided dg-ideal follows
from Lemma~\ref{lambdamodannhasonlyonesimple}.
By definition there is a dg-simple dg-module, and
$$\dgrad_2(\Lambda,d)=\bigcap_{(S,\delta)\textup{ dg-simple dg module}}\ann(S,\delta)$$
is an intersection over a non empty index set.
We have a canonical dg-ring homomorphism
$$\Lambda\lra
\prod_{{\textup{\begin{minipage}{3.2cm}$(S,\delta)$ representative of \\ an  isoclass of dg-simple\\ dg-module up to shift
\end{minipage}}}}
\Lambda/\ann(S,\delta).$$
Composing this map with the projection onto $\Lambda/\ann(S,\delta)$ is the natural projection
$\Lambda\lra\Lambda/\ann(S,\delta)$, which clearly is surjective. Further, the kernel is
precisely $\dgrad_2(\Lambda,d)$. Hence, we obtain that
$(\Lambda/\dgrad_2(\Lambda,d),\ol d)$ is a subdirect product of dg-primitive dg-algebras.

Let $(I,d)$ be a twosided dg-ideal such that $(\Lambda/I,\ol d)$ is a subdirect product
of dg-primitive dg-algebras. Then the projection onto each of the direct factors is a
dg-primitive dg-algebra quotient, yields hence a dg-simple faithful dg-module. Hence, its annihilator
contains $I$.
As $I$ is a subdirect product of all these dg-primitive algebras, it is contained in
the intersection of all the annihilators of the corresponding dg-simples.
Hence, $\dgrad_2(\Lambda,d)\subseteq I$. \dickebox

\begin{Rem}
Recall that an artinian primitive algebra is simple.
One might ask if similar properties hold in the differential graded case. One major
difficulty is that we do not have a Wedderburn Artin theorem in the dg-case.
Our favorite example $(Mat_{2\times 2}(K,d_1)$ provides a striking
observation. The dg-algebra has only one dg-simple ideal, namely the right
hand column. This is a faithful dg-simple dg-module, and hence the algebra is dg-primitive.
Since it is simple, it is also dg-simple (and therefore also dg-primitive, giving a second argument).
However, the left regular dg-module is not a direct sum of dg-simples. Actually, the
dg-left module structure is one of the type studied in Remark~\ref{exampleofdgsimpleoverdgsimple}.
Nevertheless, the algebra is acyclic.
%
\end{Rem}

}

\section{Differential Graded Orders}

\label{diffgradedordersection}

Let $R$ be a Dedekind domain with field of fractions $K$. Recall that $K$ is a flat $R$-module.

\begin{Def}
A {\em differential graded $R$-order (dg-$R$-order)} is a differential graded $R$-algebra
$(\Lambda,d)$ such that $\Lambda$ is {\red finitely generated $R$-projective} and such that
$K\otimes_RA$ is semisimple artinian as an algebra. A differential graded $R$-order is a {\em proper
differential graded $R$-order} if in addition
$H(\Lambda)$ is {\red finitely generated $R$-projective,}
\end{Def}

\begin{Def}
Let $(\Lambda,d)$ be a dg-$R$-order. A {\em differential graded $\Lambda$-lattice
(dg-$R$-lattice, or dg-lattice)} is a differential
graded $\Lambda$-module $(L,d_L)$ such that $L$ is $R$-{\red projective}. A dg-lattice $(L,d_L)$
is a {\em proper dg-lattice} if in addition $H(L)$ is $R$-{\red projective}.
\end{Def}

\begin{Rem}
Note that if $(\Lambda,d)$ is a dg-order, then a dg-lattice $(L,d_L)$ can be proper
or not. If $(\Lambda,d)$ is a proper dg-order, then a dg-lattice $(L,d_L)$ can
be proper or not.
\end{Rem}

\begin{Rem} Let $R$ be a Dedekind domain with field of fractions $K$.
Following Theorem~\ref{GarciaRozastheorem} we may as well consider
differential graded $K$-algebras $(A,d)$ which are {\red acyclic} and $\ker(d)$ semisimple.
Then we may define an $R$-order as a differential graded $R$-subalgebra $(\Lambda,d)$
projective as an $R$-module with $K\otimes_R\Lambda=A$.
This will give a different definition, implying that the homology is necessarily
$R$-torsion. Further, the fact that the algebra $A$ is very general
does not seem to allow a rich theory. Moreover, as a classical order is not a dg-order
with differential $0$ in this case, this alternative theory will not generalize
the classical theory of lattices over orders, and in particular it is unclear how
one might recover in a profitable way complexes of lattices over a classical order in any way
as differential graded module over a dg-order in this alternative definition.
\end{Rem}

\begin{Example}
We recall some examples of classical orders and explain differential graded structures on them.
\begin{enumerate}
\item For $R$ a complete discrete valuation ring
containing $\widehat\Z_3$ the $3$-adic integers, following~\cite[Example 1.2.43.4]{grouprepbuch}
the group ring $R\Sn_3$ of the
symmetric group of order $6$ is isomorphic to
$$
\{(d_1,\left(\begin{array}{cc}a_2&b_2\\c_2&d_2\end{array}\right),a_3)\in R\times Mat_{2\times 2}(R)\times R\;|\;
a_2-d_1\in 3R, c_2\in 3R, a_3-d_2\in 3R\}.
$$
This algebra is differential graded with the grading given by $deg(a_i)=deg(d_i)=0$,
$deg(c_2)=-1$ and $deg(b_2)=1$ for all $i$, using the designation of the variables as in
the above definition. The differential is given as in Example~\ref{Examplediffgradedalgebra}
by $d_x$ for any fixed $x\in R$. More precisely
$$d_x((d_1,\left(\begin{array}{cc}a_2&b_2\\c_2&d_2\end{array}\right),a_3))=
(0,\left(\begin{array}{cc}xc_2&x(d_2-a_2)\\0&xc_2\end{array}\right),0).$$
We can even consider suborders, such as those given by asking in addition that $a_2-d_2\in 3R$.

\item
Fix a prime $p>2$. {\red By a result due to Roggenkamp~\cite{GreenOrders}}
we know that the principal block of $R\Sn_p$, for $R$ being a complete discrete valuation
ring containing $\widehat \Z_p$, the $p$-adic integers and $\Sn_p$ the symmetric group of degree $p$,
is Morita equivalent to
$$
\{(d_1,\prod_{i=2}^{p-1}\left(\begin{array}{cc}a_i&b_i\\c_i&d_i\end{array}\right),a_p)
\in R\times \prod_{i=2}^{p-1}Mat_{2\times 2}(R)\times R\;|\;
a_{j+1}-d_j\in pR, c_j\in pR, \forall j\}.
$$
Again, following Example~\ref{Examplediffgradedalgebra},
for every of the $2\times 2$ matrix algebras in the product we may impose a
differential graded algebra structure with differential $d_{x_i}$ by
choosing $x_2,\dots,x_{p-1}\in R$.

\item
Fix a prime $p>2$ and let $R$ be a complete discrete valuation
ring containing $\widehat Z_p$, the $p$-adic integers. By a result due to K\"onig~\cite{Koenig}
the principal block of the Schur algebra $S_R(p,p)=End_{R\Sn_p}((R^p)^{\otimes p})$ with
parameters $(p,p)$ is known to be Morita equivalent to an algebra of similar shape.
$$
\{(\prod_{i=1}^{p-1}\left(\begin{array}{cc}a_i&b_i\\c_i&d_i\end{array}\right),a_p)
\in \prod_{i=1}^{p-1}Mat_{2\times 2}(R)\times R\;|\;
a_{j+1}-d_j\in pR, c_j\in pR, \forall j\}.
$$
The same statement on differential graded structures holds.

\item
Fix a prime $p$. By \cite{polydrozd} (\cite{drozd1,drozd2} for the case $p=2$ and $p=3$)
the category ${\mathcal F}^p_{\Z}$ of polynomial functors of degree at most $p$
from finitely generated
free abelian groups to finitely generated $\widehat \Z_p$-modules is equivalent to a module category,
having a 'principal block' and a number of trivial direct factors.
 The principal block is Morita equivalent to
$$
\{(d_0,\prod_{i=1}^{p-1}\left(\begin{array}{cc}a_i&b_i\\c_i&d_i\end{array}\right),a_p)
\in \widehat\Z_p\times \prod_{i=1}^{p-1}Mat_{2\times 2}(\widehat\Z_p)\times\widehat\Z_p\;|\;
a_{j+1}-d_j\in p\widehat\Z_p, c_j\in p\widehat\Z_p, \forall j\}.
$$
Again, following Example~\ref{Examplediffgradedalgebra},
for every of the $2\times 2$ matrix algebras in the product we may impose a
differential graded algebra structure with differential $d_{x_i}$ by
choosing $x_1,\dots,x_{p-1}\in p\widehat\Z_p$.

\end{enumerate}
\end{Example}

\begin{Rem}
Le Bruyn et al studied in \cite[Section II.4]{GradedOrders} graded orders
over Krull domains in a very general setting. They define for a Krull domain
$R$ with field of fractions $K$ and a central simple $K$-algebra $A$ a subring
$\Lambda$ with $R\subset \Lambda\subset A$ an $R$-order if each element $a$
of $\Lambda$ is integral over $R$, and if $K\Lambda=A$. If $a$ is integral
over $R$, then for the reduced trace $Tr:A\ra K$ one has $Tr(a)\in R$ and
$Tr$ induces an isomorphism $A\ra Hom_K(A,K)$ by
$A\ni a\mapsto(A\ni b\mapsto Tr(ab)\in K)\in Hom_K(A,K)$.
\end{Rem}

Recall in this context the classical result

\begin{Theorem} (cf. e.g. \cite{ReinerMO}, \cite[Theorem 8.3.7]{grouprepbuch})
\label{integralelementsformorders}
Let $R$ be a Dedekind domain of characteristic $0$ with field of fractions $K$
and let $\Lambda$ be an $R$-subalgebra
of the semisimple $K$-algebra $A$, containing a $K$-basis of $A$.
Then $\Lambda$ is an $R$-order if and only if every element of
$\Lambda$ is integral over $R$.
\end{Theorem}

Hence the definition of \cite[Section II.4]{GradedOrders} coincides with the classical
definition of an $R$-order in case of a Dedekind domain $R$.


Recall that a dg-module $(V,\delta)$ with $\delta(V_i)\subseteq V_{i+1}$ for all $i$
is {\em right bounded } if
there is $n_0\in\N$ such that $V_n=0$ for all $n>n_0$. Analogously we define {\em left bounded} dg-modules.
A dg-module is {\em bounded} if it is at once left and right bounded. Note that for a differential graded
$R$-order $(\Lambda,d)$ in a finite dimensional semisimple differential
graded $K$-algebra $(A,d)$ the lattice
$(\Lambda,d)$ is always bounded. Indeed, since $(A,d)$ is finite dimensional, $(A,d)$ is bounded,
and hence so is $(\Lambda,d)$.

\begin{Lemma} \label{existenceofdglatticesindgmodules}
Let $R$ be a Dedekind domain with field of fractions $K$.
Let $(\Lambda,d)$ be a differential graded $R$-order in a finite dimensional
differential graded  $K$-algebra $(A,d)$
and let $(V,\delta)$ be finite dimensional differential graded (hence bounded) $(A,d)$-module. Then there
is a differential graded $(\Lambda,d)$-lattice $(L,\delta)$ in $(V,d)$ such that $K\otimes_RL=V$.
\end{Lemma}

Proof. Recall that classically for an $R$-order $\Lambda$ and a $K\Lambda$-module $V$ there
is always a
full $\Lambda$-lattice $L$ in $V$. Consider $\Lambda^+:=\bigoplus_{k\geq 1}\Lambda_k$ and
$\Lambda_0^+:=\Lambda_0\oplus\Lambda^+$. Now, $\Lambda_0^+$ is a
(differential graded) subalgebra of $\Lambda$ and $\Lambda^+$ is a dg-ideal of $\Lambda_0^+$.

Let us perform the differential graded construction.
We construct $(L,\delta)$ by downward induction on the degree. Let $V_k=0\neq V_{n-1}$ for $k\geq n$
and take any
$\Lambda_0^+$-lattice in $V_{n-1}$. Its existence is a well-known classical property for
lattices.  Note that $\Lambda^+$ acts as $0$.
Hence this is automatically a $(\Lambda_0^+,d)$-lattice since $\delta(V_{n-1})=0$ and hence the
Leibniz formula automatically holds.

Let $\widetilde L_{n-2}:=\delta^{-1}(L_{n-1})\subseteq V_{n-2}$ and choose a full
$\Lambda_0$-lattice $\widehat L_{n-2}$
in $\widetilde L_{n-2}$. Then choose an $r\in R\setminus\{0\}$ such that
$r\cdot \Lambda^+\cdot \widehat L_{n-2}\subseteq L_{n-1}$.
This is possible since $\Lambda^+\cdot \widehat L_{n-2}$
is a $\Lambda^+_0$-lattice. Put $L_{n-2}:=r\cdot \widehat L_{n-2}$.
Then by construction $\delta(L_{n-2})\subseteq L_{n-1}$ and
$L_{n-2}\oplus L_{n-1}$ is a $\Lambda^+_0$-lattice,
and moreover the Leibniz formula holds since it holds for the $(A,d)$-action on $(V,\delta)$.
We suppose having constructed $L_k,L_{k+1},\dots,L_n$.
Then let $\delta^{-1}(L_k)=:\widetilde L_{k-1}$ and choose a $\Lambda_0$-lattice $\widehat L_{k-1}$
in $\widetilde L_{k-1}$.
Then there is again an $r\in R\setminus\{0\}$ such that
$$r\cdot\Lambda^+\cdot \widehat L_{k-1}\subseteq\bigoplus_{\ell=k}^nL_\ell.$$
Put $L_{k-1}:= r\cdot\widehat L_{k-1}$. Again, $\bigoplus_{\ell=k-1}^nL_\ell$
this is a full differential graded $\Lambda^+_0$-lattice. Since $(V,\delta)$ is bounded,
after a finite number of steps we constructed a Noetherian differential graded
$\Lambda^+_0$-lattice $(L,\delta)$ as $\bigoplus_{k\in\Z}L_k$.

Now, since $\Lambda$ is Noetherian as well, $\Lambda\cdot L$ is Noetherian again,
and since it is a submodule of $V$, it is $R$-torsion free. Hence $\Lambda\cdot L$ is a lattice
in $(V,\delta)$. Moreover, the Leibniz formula holds.
Since $$\delta(\lambda\cdot x)=d(\lambda)\cdot x\pm \lambda \cdot\delta(x)\in\Lambda\cdot L$$
for homogeneous elements $\lambda\in\Lambda$ and $x\in L$. This shows that
$$\delta(\Lambda\cdot L)\subseteq\Lambda\cdot L$$
is a full differential graded $(\Lambda,d)$-lattice in $(V,d)$. \dickebox

\begin{Rem}
Note that we do not really need $\Lambda_0^+$. It is possible to
first find a full and finitely generated $R$-submodule $\check L$ of $V$,
stable under $\delta$,  and then consider $\Lambda\cdot\check L$. Noetherianity and
the last argument of the above proof then
provides the result.  However, the above proof of Lemma~\ref{existenceofdglatticesindgmodules}
gives a more direct construction for coconnective dg-algebras. A similar construction
can be given for connective dg-algebras.
\end{Rem}

\medskip

In the classical theory of orders the following result is a main tool. We shall need to
transpose it to the differential graded situation.

\begin{Prop}\label{maxheredlocal}
Let $R$ be a Dedekind domain with field of fractions $K$. For each $\wp\in \textup{Spec}(R)$
and any $R$-module $M$ denote $M_\wp:=R_\wp\otimes_RM$ and $\textup{id}_{R_\wp}\otimes d:=d_\wp$ for
any map $d$, in particular the differential.
\begin{enumerate}
\item \label{local1}
If $(\Lambda,d_\Lambda)$ is a (resp. proper) dg-$R$-order in the semisimple dg-$K$-algebra $(A,d_A)$,
then $(\Lambda_\wp,1_{R_\wp}\otimes d_\Lambda)$ is a (resp. proper) dg-$R_\wp$-order
in the semisimple dg-$K$-algebra $(A,d_A)$.
\item \label{local2}
If $(L,d_L)$ is a (resp. proper) dg-$R$-lattice, then $(L_\wp,(d_L)_\wp)$ is a (resp. proper)
dg-$(\Lambda_\wp,(d_\Lambda)_\wp)$-lattice.
\item \label{local3}
For any (resp. proper) dg-$(\Lambda,d_\Lambda)$-lattice $(L,d_L)$ we have
$$L=\bigcap_{\wp\in \textup{Spec}(R)}(L_\wp,(d_L)_\wp)$$
where the intersection is taken inside $K\otimes_RL$.
\item \label{local4}
Fix a dg-$(A,d_A)$-module $(V,d_V)$.
For each $\wp\in\textup{Spec}(R)$ fix (resp. proper)
$(\Lambda_\wp,(d_\Lambda)_\wp)$-lattices $(M(\wp),(d_M(\wp)))$ such that
$(KM(\wp),(Kd_M(\wp)))=(V,d_V)$ for all $\wp\in\textup{Spec}(R)$. Suppose moreover
that there is a (resp. proper) dg-$(\Lambda,d_\Lambda)$-lattice $(N,d_N)$
such that $(N_\wp,(d_N)_\wp)=(M(\wp),d_M(\wp))$ for all but a finite number of
$\wp\in\textup{Spec}(R)$.
Then there is a (resp. proper) dg-$(\Lambda,d_\Lambda)$-lattice $(L,d_L)$ with
$(L_\wp,(d_L)_\wp)=(M(\wp),d_M(\wp))$ for all $\wp\in\textup{Spec}(R)$.
\item \label{local5}
The analogous statements hold replacing the localisation by the completion.
\end{enumerate}
\end{Prop}

Proof.
Item (\ref{local1}) follows from the classical non dg-statement
(cf e.g. \cite[Proposition 8.1.14]{grouprepbuch}).

For item (\ref{local2}) we need to see that $H(L_\wp,(d_L)_\wp)$ is $R_\wp$-torsion free
if $H(L,d_L)$ is $R$-torsion free.
But this follows from the fact that localisation is flat (cf e.g. \cite[Lemma 6.5.7]{grouprepbuch}).

Item (\ref{local3}) is again a direct consequence of the classical non dg-statement
(cf e.g. \cite[Proposition 8.1.14]{grouprepbuch}).

As for item (\ref{local4}) we first get a lattice $L$ again from the
classical situation
(cf e.g. \cite[Proposition 8.1.14]{grouprepbuch}). The differential is fixed as
the restriction of the differential $d_V$ on $L$.
We need to verify that $d_V(L)\subseteq L$. But this is true at every prime, i.e.
$$d_V(L_\wp)=d_M(\wp)(M(\wp))\subseteq M(\wp)=L_\wp.$$
Hence by the non dg-version of item (\ref{local3}) we have $d_V(L)\subseteq L$.

As for item (\ref{local5}) we first see that the non dg-version is again
classical (cf e.g. \cite[Proposition 8.1.14]{grouprepbuch}).
The analogous of the second statement first uses the localisation, then
the uniqueness and existence of the differential follows by continuity.
The third item is clear by the localisation case.
As for the fourth item we first restrict from the completion to $A$. Then,
in a second step we use continuity again.
\dickebox

\medskip

An important tool in case of lattices for orders is the conductor. Recall that given a
Dedekind domain $R$ with field of fractions $K$ and a finite dimensional semisimple
$K$-algebra $A$, then for an $R$-lattice $L$ in $A$ with $KL=A$ we define
$${\mathcal O}_\ell(L):=\{\lambda\in A\;|\;\lambda L\subseteq L\}
\text{ and }{\mathcal O}_r(L):=\{\lambda\in A\;|\; L\lambda\subseteq L\}.$$
Then by \cite[Lemma 8.3.14]{grouprepbuch} we get that ${\mathcal O}_\ell(L)$
and ${\mathcal O}_r(L)$ are $R$-orders
in $A$.

We can show a dg-version of this lemma.

\begin{Lemma}\label{leftonductorlemmadg}
Let $R$ be a Dedekind domain with field of fractions $K$ and {\red let $(A,d)$ be a
differential graded $K$-algebra $(A,d)$, semisimple as an algebra.}
Let $(L,d_L)$ be a differential graded $R$-lattice in $(A,d)$  such that $K\cdot (L,d_L)=(A,d)$.
Then $({\mathcal O}_\ell(L),d|_{{\mathcal O}_\ell(L)})$ is a differential graded $R$-order in $(A,d)$.
Similar statements hold for ${\mathcal O}_r(L)$.
\end{Lemma}

Proof. By symmetry it is enough to consider the case $({\mathcal O}_\ell(L),d|_{{\mathcal O}_\ell(L)})$.
By \cite[Lemma 8.3.14]{grouprepbuch} we see that ${\mathcal O}_\ell(L)$ is an $R$-order in $A$.
{\red If $L$ is graded, then also ${\mathcal O}_\ell(L)$ is graded. }
We need to see that it is differential graded.
Let $\lambda\in {\mathcal O}_\ell(L)$ be homogeneous. Then $\lambda L\subseteq L$.
We need to see that $d(\lambda) L\subseteq L$.
Since $L\subseteq A$, and since $(A,d)$ is a differential graded algebra, we
have for any $x\in L$ and homogeneous $\lambda\in {\mathcal O}_\ell(L)$
$$d(\lambda x)=d(\lambda)x+(-1)^{|\lambda|}\lambda d(x)$$
Hence $$d(\lambda)x=d(\lambda x)-(-1)^{|\lambda|}\lambda d(x)$$
Since $\lambda\in {\mathcal O}_\ell(L)$ and since $d(L)\subseteq L$, using that $(L,d)$ is
a dg-lattice, we have $\lambda d(x)\in L$. Since $\lambda\in {\mathcal O}_\ell(L)$, and since $x\in L$,
we get $\lambda x\in L$, and hence $d(\lambda x)\in L$ again. The differential on ${\mathcal O}_\ell(L)$
is the restriction of the differential $d$ of $A$ to  ${\mathcal O}_\ell(L)$, and hence
the defining equation on products holds.
 Hence $({\mathcal O}_\ell(L),d|_{{\mathcal O}_\ell(L)})$
is a dg-order in $(A,d)$. \dickebox

\begin{Rem} \label{leftorder-dg-version-cannot-be-proper}
{\red If the category of dg-modules over $(A,d)$ is semisimple,} then by \cite{Tempest-Garcia-Rochas}
we have that $H(A,d)=0$ and therefore any order in $(A,d)$ is either {\red acyclic} or
has $R$-torsion homology.
In particular, in this case, if the conductor is not {\red acyclic}, then
 $({\mathcal O}_\ell(L),d|_{{\mathcal O}_\ell(L)})$
cannot be a proper differential graded $R$-order in $(A,d)$, whatever the choice of
a dg-lattice $L$ may be.
\end{Rem}

\section{On maximal differential graded orders}

\label{maximaldiffgradedorders}

Let $R$ be a Dedekind domain with field of fractions $K$.
It is a classical fact that for any semisimple $K$-algebra $A$ and an $R$-order $\Lambda$
in $A$ there is a maximal $R$-order $\Gamma$ containing $\Lambda$. Maximal orders
have many striking properties, and behave very much as the base ring $R$ (cf \cite{ReinerMO}).
We shall consider maximal dg-orders.

\begin{Theorem}\label{gradedmaximalordersexist}
Let $R$ be a Dedekind domain with field of fractions $K$ {\red of characteristic $0$}
and let $(\Lambda,d)$ be a
differential graded $R$-order in the semisimple finite dimensional
differential graded $K$-algebra $(A,d)$.
Then there is a differential graded  $R$-order $(\Gamma,d)$, {\red which is maximal
with respect to being a dg-order and} containing $(\Lambda,d)$.
If $(\Lambda,d)$ is a proper dg-order,  then there is a proper
differential graded  $R$-order $(\Gamma_p,d)$ {\red which is maximal
with respect to being a proper dg-order and} containing $(\Lambda,d)$.
\end{Theorem}

Proof. The set of dg-orders in $(A,d)$ is partially ordered.
We need to make precise the partial ordering we consider.
Let $(\Lambda_1,d_1)$ and $(\Lambda_2,d_2)$ be dg-$R$-orders in a semisimple
dg-$K$-algebra $(A,d)$. Then
$(\Lambda_1,d_1)\leq (\Lambda_2,d_2)$ if $\Lambda_1\subseteq\Lambda_2$.
Since $d_1=d|_{\Lambda_1}$ and $d_2=d|_{\Lambda_2}$, we get that then automatically
$d_1=d_2|_{\Lambda_1}$.

Consider
$${\mathcal X}:=\{(\Sigma,d)\;|\;(\Lambda,d)\leq(\Sigma,d)\;\textup{ and }
(\Sigma,d) \textup{ is a dg-$R$-order in $(A,d)$}\}.$$
Since $(\Lambda,d)\in{\mathcal X}$, this set ${\mathcal X}$ is not empty.
Let ${\mathcal Y}$ be a totally ordered subset in ${\mathcal X}$ and put
 $\Gamma:=\bigcup_{(\Delta,d)\in{\mathcal Y}}\Delta$. It is graded since
the embedding of one order in the other in ${\mathcal Y}$ preserves the grading, and hence so is
the union. It allows a differential by
$d|_\Gamma$ and if $\gamma\in\Gamma$, then $\gamma\in\Delta$ for some
$\Delta\in{\mathcal Y}$, and therefore $d(\gamma)\in\Delta\subseteq\Gamma$.
Moreover, since $\gamma\in\Delta$, and since $\Delta$ is an order,
Theorem~\ref{integralelementsformorders} shows that
$\gamma$ is integral, and, using Theorem~\ref{integralelementsformorders} again,
$\Gamma$ is an order in $A$.
Hence $(\Gamma,d)\in{\mathcal X}$ and
$(\Gamma,d)$ dominates any element in $\mathcal Y$. By Zorn's lemma there are maximal elements in
$\mathcal X$. Any such element is a maximal differential graded order containing $(\Lambda,d)$.

Now consider
$$
\widehat{\mathcal X}:=\{(\Sigma,d)\;|\;(\Lambda,d)\leq(\Sigma,d)\;\textup{ and }
(\Sigma,d) \textup{ is a proper dg-$R$-order in $(A,d)$ and $\Sigma$}\}$$
and if $(\Lambda,d)$ is
a proper dg-$R$-order in $(A,d)$, then  $(\Lambda,d)\in\widehat{\mathcal X}$. Hence
this set is not empty neither.  The first steps are as above.
Let $\widehat{\mathcal Y}$ be a totally ordered subset in $\widehat{\mathcal X}$ and put
$\Gamma:=\bigcup_{(\Delta,d)\in\widehat{\mathcal Y}}\Delta$. It is graded since
the embedding of one order in the other in ${\mathcal Y}$ preserves the grading,
and hence so is the union. It allows a differential by
$d|_\Gamma$ and if $\gamma\in\Gamma$, then $\gamma\in\Delta$ for some $\Delta\in\widehat{\mathcal Y}$.
Hence $\gamma$ is integral, and therefore $\Gamma$ is an order in $A$.
Since $R$ is Noetherian (cf \cite[Lemma 7.5.3]{grouprepbuch}), $\Gamma$ is Noetherian as well, and therefore
$\Gamma$ actually in $\widehat{\mathcal Y}$. Therefore $(\Gamma,d)$ is a proper dg-order.
By Zorn's lemma, $\widehat{\mathcal X}$ contains maximal elements.
This proves the statement.
\dickebox

%
%
%

\medskip

{\red We call an $R$-order, which is maximal with respect to being a dg-order in a fixed
dg-algebra, a dg-maximal dg-order.}

\begin{Cor}\label{maximalordersgivedifferentialgradedorders}
Let $R$ be a Dedekind domain with field of fractions $K$, let $(A,d)$ be a finite dimensional
semisimple differential graded $K$-algebra, and let $\Lambda$ be a maximal $R$-order in $A$
and suppose that $d(\Lambda)\subseteq\Lambda$.
Then $(\Lambda,d|_\Lambda)$ is a {\red dg-maximal} differential graded $R$-order in $(A,d)$.
\end{Cor}

Proof. The hypotheses imply that $(\Lambda,d|_\Lambda)$ is a differential $R$-graded order.
If $(\Gamma,d|_\Gamma)$ is a differential graded order containing $(\Lambda,d|_\Lambda)$,
then $\Gamma$ is an $R$-order containing $\Lambda$.
Hence, since $\Lambda$ is a {\red dg-}maximal order, $\Lambda=\Gamma$ and therefore
$(\Lambda,d|_\Lambda)$ is a {\red dg-}maximal differential $R$-graded order.
This shows the statement. \dickebox

\medskip

\begin{Rem}\label{not-proper-dg-order}
Note that in the situation of Corollary~\ref{maximalordersgivedifferentialgradedorders}
there is no reason why $H(\Lambda,d_\Lambda)$ should be $R$-projective, and indeed
$H(A,d)=0$, as in Example~\ref{Examplediffgradedalgebra}, is possible and implies this case.
Recall from Remark~\ref{tempest-GarciaRozas}
that orders in algebras {\red with
semisimple category of dg-module} are actually necessarily of this form.
Let $(\Lambda,d)$ be an order in $(A,d)$ with $H(\Lambda,d)=0$. Then
Aldrich and Garcia-Rocas show in \cite[Theorem 4.7]{Tempest-Garcia-Rochas} that
in this case the category of differential graded $(\Lambda,d)$-modules is equivalent with the
category of graded $\ker(d)$-modules.
\end{Rem}

\begin{Example}
Let $R$ be an integral domain with field of fractions $K$. We consider the differential graded
semisimple $K$-algebra $(Mat_{2\times 2}(K),d_x)$ from Example~\ref{Examplediffgradedalgebra}.
Recall that for each $x\in K$ there is a differential $d_x$ on $Mat_{2\times 2}(K)$
with the grading chosen in Example~\ref{Examplediffgradedalgebra}.

Choose $x=1$ for the moment.
Then  $Mat_{2\times 2}(R)$ is a maximal differential graded $R$-order, and actually
a proper differential graded $R$-order since the homology of $(Mat_{2\times 2}(R),d_1)$ is $0$.

However, choosing an element $x\in R\setminus R^\times$, then again
$Mat_{2\times 2}(R)$ is a maximal differential graded $R$-order. The homology in degree $1$ is $R/xR$,
as is easily seen. The kernel of the differential in degree $0$ is
$C_0:=\{\left(\begin{array}{cc}a&0\\0&a\end{array}\right)\;|\;a\in R\}$
and the image if the differential from degree $-1$ is $xC_0$. Clearly, the differential
in degree $-1$ is injective. Hence $H_*(Mat_{2\times 2}(R),d_x)\simeq R/xR\oplus R/xR$,
where the first copy is in degree $0$ and the second copy is in degree $1$.
Therefore  $H_*(Mat_{2\times 2}(R),d_x)\simeq (R/xR)[\epsilon]/\langle \epsilon^2\rangle$
where $\epsilon$ is an element in degree $1$.

Note that $Mat_{2\times 2}(R)$ is hereditary, as an order, whereas, as soon
as $x\not\in R^\times$, the homology algebra $H_*(Mat_{2\times 2}(R),d_x)$
is not, even of infinite global dimension.
%
\end{Example}

\begin{Example} Let $R$ be an integral domain with field of fractions $K$.
We consider again the differential graded
semisimple $K$-algebra $(Mat_{2\times 2}(K),d_x)$ from Example~\ref{Examplediffgradedalgebra}.
If $R=\Z$ and $x=\frac 12$, then for $\Lambda=Mat_{2\times 2}(R)$ we obtain
$d_x(\Lambda)$ is not a subset of $\Lambda$. Actually,
$$d_{\frac{1}{2}}(Mat_{2\times 2}(\Z))=\frac12\Z\left(\begin{array}{cc}1&0\\ 0&1
\end{array}\right)+\frac12\Z\left(\begin{array}{cc}0&1\\ 0&0
\end{array}\right).$$
Since
$$\left(\begin{array}{cc}\frac12&0\\0&\frac12\end{array}\right)\cdot \left(\begin{array}{cc}0&\frac12\\0&0\end{array}\right)=
\left(\begin{array}{cc}0&\frac14\\0&0\end{array}\right)$$
and iterating, we have that for each integer $n$ the element
$\left(\begin{array}{cc}0&\frac1{2^n}\\0&0\end{array}\right)$ is in the ring generated by $d_{\frac12}(Mat_{2\times 2}(\Z))$.
Hence, there is no differential graded order $(\Gamma,d|_\Gamma)$ such that
$Mat_{2\times 2}(K)\subseteq\Gamma$.
However, for any $x\in K\setminus\{0\}$ the set
$$
\Lambda:= \left(\begin{array}{cc}
R&xR\\\frac{1}{x}R&R
\end{array}\right)
$$
is a subring of $Mat_{2\times 2}(K)$ and is stable under $d_x$. Further,
if $R$ is a Dedekind domain, then $xR$ is projective, and hence $\Lambda$ is an
$R$-order if $R$ is a Dedekind domain.
We observe that
$$\left(\begin{array}{cc}\frac1x &0\\0&1\end{array}\right)\cdot\Lambda\cdot
\left(\begin{array}{cc}x &0\\0&1\end{array}\right)=Mat_{2\times 2}(R)$$
and hence $\Lambda$ is conjugate to the maximal order $Mat_{2\times 2}(R)$. Therefore,
using Corollary~\ref{maximalordersgivedifferentialgradedorders}, $(\Lambda,d_x)$
is a maximal differential graded order.
We get $H(\Lambda,d_x)=0$, and hence
this is actually a proper differential graded $R$-order.
\end{Example}

\begin{Question}
If $(\Lambda,d)$ is a {\red dg-}maximal differential graded order, can one show that $\Lambda$
is a maximal order.
\end{Question}

\begin{Question}
If $(\Lambda,d)$ is a {\red dg-}maximal proper differential graded order, can we show that
$H(\Lambda,d)$ is hereditary?
\end{Question}

\section{Class groups of differential graded orders}

\label{classgroupsection}


Recall that a differential graded $R$-order $(\Lambda,d)$
in a semisimple differential graded $K$-algebra $(A,d)$ is at first an
$R$-order $\Lambda$ in a semisimple $K$-algebra $A$. The locally free class group
of an order proved to be a useful invariant in integral representation theory. We
want to provide a definition and first properties of a locally free class group
of differential graded orders.

We first need to elaborate on what we mean by a free differential graded module.
Let $(\Lambda,d)$ be a differential graded algebra.
A differential graded $(\Lambda,d)$-module $(L,\delta)$ is free of rank $n\in\N$ if
$(L,\delta)\simeq\bigoplus_{i=1}^n(\Lambda,d){\red [k_i]}$ {\red for some integers $k_i$},
as differential graded $(\Lambda,d)$-modules.
{\red Later we shall mainly consider the case $k_i=0$ for all $i$, and say that
such a dg-module is degree $0$-free of rank $n$.}

Since $d(1)=0$ in $\Lambda$, we get that for an
{\red homomorphism $\varphi:\bigoplus_{i=1}^n(\Lambda,d)[k_i]\lra (L,\delta)$}
we have $\varphi(0,\dots,0,1,0,\dots,0)\in\ker(\delta)$, where $1$ is in position $i$,
for each position $i$. However, for any {\red homogeneous}
$z\in\ker(\delta)$ we get that $\varphi(\lambda):=\lambda z$
defines a homomorphism $(\Lambda,d){\red [-|z|]}\lra (L,\delta)$.

\begin{Lemma} \label{freeintermsofunits}
Let $R$ be a Dedekind domain with field of fractions $K$.
Let $(\Lambda,d)$ be a differential graded algebra, and
denote $(A,d):=(K\otimes_R\Lambda,\textup{id}_K\otimes d)$.
Then the set of free differential graded $(\Lambda,d)$-ideals
is in bijection with the group of left
regular {\red homoegenous} elements in $\ker(d)\subseteq A$
modulo the subgroup $(\ker(d)\cap\Lambda^\times)$.
\end{Lemma}

Proof.
We now consider dg-ideals
$(I,d)$ of $(\Lambda,d)$. Such an ideal is free if there is an isomorphism $\varphi:(\Lambda,d){\red [k]}\lra (I,d)$
and this is equivalent with the choice of a $z\in (I\cap\ker(d))$ {\red homogeneous
of degree $k$} such that $I=\Lambda z$ and such that
$\lambda z=0\Rightarrow \lambda=0$. We hence obtain that any left regular {\red homogeneous}
generator of $I$ in the cycles
gives an isomorphism, and hence there is a surjective map from left regular {\red homogeneous}
elements in the cycles to the
set of rank one free dg $(\Lambda,d)$ ideals.
Two such  left regular elements {\red homogeneous} $z_1,z_2$ of $I$ in the cycles give the same ideal if and only if
$\Lambda z_1=\Lambda z_2$. This is equivalent with $z_1=\lambda_2z_2$ and $z_2=\lambda_1z_1$
for {\red homogeneous} units $\lambda_1,\lambda_2\in\Lambda$. However, $z_1,z_2\in\ker(d)$ shows that
$$0=d(\lambda_2z_2)=d(\lambda_2)z_2\pm\lambda_2d(z_2)=d(\lambda_2)z_2$$
and the fact that $z_2$ is regular gives that $\lambda_2\in\ker(d)$. Likewise $\lambda_1\in\ker(d)$.

Conversely, if $\lambda\in\ker(d)$ {\red is a homogeneous unit}, then $\Lambda z=\Lambda\lambda z$ as
differential graded ideal. Indeed, if $\lambda\in\ker(d)$ and $z\in\ker(d)$, then
$$
d(\mu \lambda z)=
d(\mu)\lambda z+(-1)^{|\mu|}\mu d(\lambda)z+(-1)^{|\mu{\red |+|}\lambda|}\mu\lambda d(z)
=d(\mu)\lambda z
$$
and hence mapping $\mu z$ to $\mu\lambda z$ for any $\mu\in\Lambda$ is an isomorphism
of differential graded ideals. \dickebox

\bigskip

We shall need to study invertible elements in the subring of cycles.

\begin{Lemma} \label{unitsincyclesareunitsofcycles}
Let $(\Lambda,d)$ be a differential graded algebra. Then $\ker(d)\cap\Lambda^\times=\ker(d)^\times$.
\end{Lemma}

Proof.
We have seen in Corollary~\ref{cyclesaresubalgebras} that $\ker(d)$ is a graded subalgebra.
Trivially  $\ker(d)\cap\Lambda^\times\supseteq\ker(d)^\times$.
Further, if $u$ is {\red homogeneous and }invertible in $\Lambda$, and if $u\in\ker(d)$, then also $u^{-1}\in\ker(d)$.
Indeed, $0=d(1)=d(uu^{-1})=d(u)u^{-1}\pm ud(u^{-1})$ and hence $d(u^{-1})=0$ since $u$ is
invertible. Hence  $\ker(d)\cap\Lambda^\times\subseteq\ker(d)^\times$. \dickebox

\bigskip

Let $R$ be a Dedekind domain with field of fractions $K$.
Suppose now that $(A,d)$ is a  differential graded $K$-algebra, {\red which is assumed to be
left and right artinian as an algebra}.
By Wedderburn's theorem, in {\red an artinian} semisimple $K$-algebra an element
$u$ is left regular if and only if it is represented by a tuple of non singular matrices, whence
an invertible element. {\red In any left and right
artinian ring a left or right regular element is invertible. Indeed, by the above, this is
clear modulo the
Jacobson radical. Then, for an artinian ring the Jacobson radical is
nilpotent, and hence any lift of a unit in the semisimple quotient to an element in the
artinian ring is again a unit.}

Let $(\Lambda,d)$ be a differential graded subalgebra of $(A,d)$
such that $\Lambda$ is $R$-projective and $K\otimes_R\Lambda=A$.
As in the classical case we need to work with
two versions, the localised version and the completed version. Let $X_\wp$ be the localisation at a
prime $\wp$, and denote by $\widehat X_\wp$ the completion at the prime $\wp$. Hence we denote
\begin{align*}
J(A,d):=\lefteqn{\{(u(\wp))_{\wp\in Spec(R)}\in \prod_{\wp\in Spec(R)}\left(A_\wp^\times\cap\ker(d)\right)\;|\;}\\
&u(\wp)\in\Lambda_\wp^\times\text{ for almost all $\wp\in Spec(R)$ {\red and   $u_\wp$ homogeneous for all $\wp$}}\}.
\end{align*}
Note that this definition does not depend on the choice of $\Lambda$.
Define its subgroup
$$U(\Lambda,d):={\red J(A,d)\cap}\prod_{\wp\in Spec(R)}(\Lambda_\wp^\times\cap\ker(d)),$$
which does depend on $\Lambda$.
Likewise
\begin{align*}
\widehat J(A,d):=\lefteqn{\{(u(\wp))_{\wp\in Spec(R)}\in
\prod_{\wp\in Spec(R)}\left(\widehat A_\wp^\times\cap\ker(d)\right)\;|\;}\\ &
u(\wp)\in\widehat\Lambda_\wp^\times\text{ for almost all $\wp\in Spec(R)$ {\red and   $u_\wp$ homogeneous for all $\wp$}}\}.
\end{align*}
and its subgroup
$$\widehat U(\Lambda,d):={\red \widehat J(A,d)\cap}\prod_{\wp\in Spec(R)}(\widehat \Lambda_\wp^\times\cap\ker(d)).$$
{\red Note that if $A$ is finite dimensional, then any homogeneous $u_\wp$ being a factor in an
element in $\widehat J(A,d)$ (resp. in $J(A,d)$)  has to be in degree $0$. }
Then the set of left classes
$$U(\Lambda,d)\backslash J(A,d)$$
is in bijection with the set of locally free differential graded fractional ideals of $(\Lambda,d)$
and
$$\widehat U(\Lambda,d)\backslash \widehat J(A,d)$$
is in bijection with the set of completion locally free differential graded
fractional ideals of $(\Lambda,d)$.

We associate to a representative $\alpha=(\alpha(\wp))_{\wp\in Spec(R)}$ of a
class $U(\Lambda,d)\alpha\in U(\Lambda,d)\backslash J(A,d)$ the fractional ideal
$${\red \Lambda\alpha:=A\cap}\bigcap_{\wp\in Spec(R)}\Lambda_\wp\cdot\alpha(\wp).$$
Likewise we get the completed version.
Then, ${\red \Lambda\alpha\simeq \Lambda\beta}$ as fractional differential graded ideal if
and only if there is {\red a homogeneous} $x\in A^\times\cap\ker(d)$  {\red of degree $0$} with $\alpha=\beta x$.
If $A$ is semisimple, then we may apply
Proposition~\ref{maxheredlocal} we have that
$${\red (\Lambda\alpha)}_\wp={\red \Lambda_\wp}\cdot\alpha(\wp)$$
for all $\wp\in Spec(R)$.

{\red We now define a category $(\Lambda,d)-{\red LF_0}-dgmod$. This is defined to be the
full subcategory of $(\Lambda,d)-dgmod$ containing all
locally degree $0$ rank one free differential graded $(\Lambda,d)$-modules, and all direct
sums of these objects.

%
}

\medskip

In the definition below we shall need to consider $K_0((\Lambda,d)-{\red LF_0}-dgmod)$.
This is defined as the quotient of the free abelian group on
isomorphism classes of objects in $(\Lambda,d)-{\red LF_0}-dgmod$
modulo the relation $[(X,\delta_X)]-[(Y,\delta_Y)]-[(Z,\delta_Z)]$ whenever
{\red $(X,\delta_X)\simeq (Y,\delta_Y)\oplus (Z,\delta_Z)$. We further define
$G_0((\Lambda,d)-{\red LF_0}-dgmod)$, which is the quotient of the free abelian group on
isomorphism classes of objects in $(\Lambda,d)-{\red LF_0}-dgmod$
modulo the relation $[(X,\delta_X)]-[(Y,\delta_Y)]-[(Z,\delta_Z)]$ whenever }
there is a
short exact sequence
$$0\lra (Y,\delta_Y)\lra (X,\delta_X)\lra (Z,\delta_Z)\lra 0$$
{\red  in $(\Lambda,d)-{\red LF_0}-dgmod$, considered as a subcategory of dg-modules.
Here, denote by $[(M,\delta)]$ the image of a locally {\red degree $0$-}free dg-module in $K_0((\Lambda,d)-LF-dgmod)$, respectively $G_0((\Lambda,d)-LF-dgmod)$.}
Note that this setting makes sense for any differential graded $R$-algebra $(\Lambda,d)$
whenever $R$ is a Dedekind domain with field of fractions $K$.

\begin{Def}
Let $R$ be a Dedekind domain with field of fractions $K$.
\begin{itemize}
\item
Let $(\Lambda,d)$ be a differential graded order in the
semisimple differential graded $K$-algebra $(A,d)$.
\begin{itemize}
\item
The group ${\red \widehat I(\Lambda,d):=}\widehat U(\Lambda,d)\backslash \widehat J(A,d)/(A^\times\cap\ker(d))$
is the group of {\em completion differential graded id\`eles }
of the dg-order $(\Lambda,d)$.
\item
The group ${\red  I(\Lambda,d):=}U(\Lambda,d)\backslash J(A,d)/(A^\times\cap\ker(d))$ is the group of {\em
differential graded id\`eles }
of the dg-order $(\Lambda,d)$.
\end{itemize}
\item If $(\Lambda,d)$ is a differential graded $R$-algebra, then
the {\em class group $Cl( (\Lambda,d))$
of  $(\Lambda,d)$} is the subgroup of $G_0((\Lambda,d)-{\red LF_0}-dgmod)$ generated by
elements $[L]-[\Lambda]$ for {\red degree $0$-}locally free differential graded
$(\Lambda,d)$-lattices $(L,\delta)$ of rank $1$.
{\red \item If $(\Lambda,d)$ is a differential graded $R$-algebra, then
the {\em id\`ele class group $Cl^{(I)}( (\Lambda,d))$
of  $(\Lambda,d)$} is the subgroup of $K_0((\Lambda,d)-{\red LF_0}-dgmod)$ generated by
elements $[L]-[\Lambda]$ for {\red degree $0$-}locally free differential graded
$(\Lambda,d)$-lattices $(L,\delta)$ of rank $1$.}
\end{itemize}
\end{Def}

As in the classical case we get

\begin{Theorem}\label{idelegroupisisotokzeroclassgroup}
Let $R$ be a Dedekind domain with field of fractions $K$. Let $(\Lambda,d)$ be a
differential graded order in the
semisimple differential graded $K$-algebra $(A,d)$.
Then {\red in  $G_0((\Lambda,d)-LF_0-dgmod)$ we have}
$$[\Lambda\alpha]+[\Lambda\beta]-2[\Lambda]=[\Lambda\alpha\beta]-[\Lambda]$$
for any two completion differential graded id\`eles $\alpha$ and $\beta$.
In particular,
{\red there is a surjective group homomorphism $\Phi$
from the group of completion differential graded id\`eles to the
differential graded class group given by $\Phi(\alpha)=[\Lambda\alpha]-[\Lambda]$.}
\end{Theorem}

Proof.
We need to verify that the constructions of id\`eles in \cite[Theorem 8.5.11]{grouprepbuch}
do not lead out of $\ker(d)$. Besides this the proof of \cite[Theorem 8.5.11]{grouprepbuch}
holds in the dg-concept
verbatim as in the non dg-concept, until the very last argument.
Let us go through the arguments of \cite[Theorem 8.5.11]{grouprepbuch}.
As in the original proof we may replace $\alpha$ by a version
after having multiplied by an $r\in R$ such that $\alpha$ is an integral id\`ele.
Again there is an integer $k$ such that $\alpha_\wp^{-1}\in\wp^{-k}\widehat\Lambda_\wp$
and an integer $t$ and $x\in A^\times\cap\ker(d)$ such that
$\beta_\wp x-1\in\wp^t\widehat\Lambda_\wp$, and again that $\beta x$ is an integral id\`ele.
We observe that $\beta x\in\ker(d)$.
Replacing $\beta$ by $\beta x$ is possible since we do not quit $\ker(d)$ by this operation.
Again, by the same computation $\alpha_\wp\beta_\wp\alpha_\wp^{-1}\beta_\wp^{-1}$
is invertible and in $\ker(d)$ as $\alpha$ and $\beta$ are.
By Lemma~\ref{unitsincyclesareunitsofcycles} units which are cycles are precisely the
units in the ring of cycles.
The morphism $f$ in \cite[page 310]{grouprepbuch} is easily seen to be a morphism
of dg-modules. The fact that $\ker(f)=\Lambda\alpha\beta$ still holds by the same proof.
The fact that $\im(f)=\Lambda$ still holds also in the dg-version by the very same proof.
Hence there is a short exact sequence of locally free dg-modules
$$0\lra\Lambda\alpha\beta\lra \Lambda\alpha\oplus\Lambda\beta\stackrel{f}\lra \Lambda\lra 0$$
and $\Lambda\alpha\beta\simeq\Lambda\beta\alpha$.
%
%
{\red All these terms are in the category
$(\Lambda,d)-{\red LF_0}-dgmod$ and we hence get a short exact sequence n the category
giving a relation in the relevant Grothendieck group.}
{\red Further, this short exact sequence shows
$$[\Lambda\alpha\beta]+[\Lambda]=[\Lambda\alpha]+[\Lambda\beta]$$
in the class group.
Recall that the group law in $Cl(\Lambda,d)$ is the group law of
the Grothendieck group, where we consider the Grothendieck group modulo direct sums.
Hence
\begin{eqnarray*}
\widehat U(\Lambda,d)\backslash \widehat J(\Lambda,d)/(A^\times\cap\ker(d))&\stackrel{\Phi}{\lra} &Cl(\Lambda,d)\\
\alpha&\mapsto&[\Lambda\alpha]-[\Lambda]
\end{eqnarray*}
is a group homomorphism. Indeed,
\begin{eqnarray*}
\Phi(\alpha\beta)&=&[\Lambda\alpha\beta]-[\Lambda]\\
&=&\left([\Lambda\alpha\beta]+[\Lambda]\right)-2[\Lambda]\\
&=&\left([\Lambda\alpha]+[\Lambda\beta]\right)-2[\Lambda]\\
&=&\left([\Lambda\alpha]-[\Lambda]\right)+\left([\Lambda\beta]-[\Lambda]\right)\\
&=&\Phi(\alpha)+\Phi(\beta)
\end{eqnarray*}.
By
construction of $Cl(\Lambda,d)$  we see that $\Phi$ is surjective.}
This proves the theorem.
\dickebox

\begin{Rem}
Note that $\Lambda$ is projective, and therefore
$$0\lra\Lambda\alpha\beta\lra \Lambda\alpha\oplus\Lambda\beta\lra \Lambda\lra 0$$
splits as $\Lambda$-modules.
However, a differential graded module $(M,\delta)$ is projective
in the category of dg-modules if and only if $(M,\delta)$ is {\red acyclic} (cf \cite{Tempest-Garcia-Rochas}).
{\red Consider the short exact sequence
$$0\lra\Lambda\gamma\stackrel\iota\lra \Lambda\alpha+\Lambda\beta\stackrel f\lra\Lambda\lra 0$$
of dg-modules over $(\Lambda,d)$. Here, $f$ is just the map sending the two components to a sum,
and $\gamma=\alpha\beta$. Consider now the completion at primes $\wp$. We only need to prove that
the sequence splits at all completions to show that the class in the $Ext$-group is $0$.

If we could have Theorem~\ref{idelegroupisisotokzeroclassgroup} for the Grothendieck group
$K_0((\Lambda,d)-{\red LF_0}-dgmod)$ instead, then we had that
the group of completion differential graded id\`eles is isomorphic to
the group $Cl^{(I)}(\Lambda,d)$ and this group actually parameterizes
stable isomorphisms of these lattices, such as in the classical situation.
}

\end{Rem}

\begin{Rem}
If $d=0$ {\red and the grading is trivial}, the concepts of class groups coincide and
$Cl(\Lambda,0)=Cl(\Lambda)$.
\end{Rem}

\begin{Cor}
Let $R$ be a Dedekind domain with field of fractions $K$ and let $(A,d)$ be a differential
graded algebra, semisimple as an algebra. Let $(\Lambda,d)$ be a differential graded order in $(A,d)$.
Then {\red $$\widehat U(\Lambda,d)\backslash {\red{\widehat J(A,d)/(A^\times\cap \ker(d)^\times)}}\simeq \widehat U(\ker(d))\backslash {\red{\widehat J(K\cdot \ker(d))/(K\cdot \ker(d)^\times)}}.$$}
\end{Cor}

Proof. 
By Lemma~\ref{unitsincyclesareunitsofcycles} we have $\widehat U(\Lambda,d)=\widehat U(\ker(d))$
and $A^\times\cap\ker(d)=\ker(d)^\times.$
\dickebox

\bigskip

More generally, we get the following

\begin{Prop} \label{ordinaryclassgroupdgclassgroup}
Let $R$ be a Dedekind domain with field of fractions $K$.
Let $(A,d)$ be a differential graded semisimple finite dimensional $K$-algebra and let
$(\Lambda,d)$ be a differential graded $R$-order in $(A,d)$.
Then there is a group homomorphism
{\red $$\widehat U(\Lambda,d)\backslash {\red{\widehat J(A,d)/(A^\times\cap \ker(d)^\times)}}\lra Cl(\Lambda).$$}
\end{Prop}

Proof.
Indeed, since $A^\times\cap\ker(d)\subseteq A^\times$ and since
$\Lambda^\times\cap\ker(d)\subseteq \Lambda^\times$, we get a group homomorphism
$J(A,d)\lra J(A)$. Further, by the same argument the same embedding
shows that $U(\Lambda,d)$ maps to $U(\Lambda)$ and hence we get a group homomorphism
$U(\Lambda,d)\backslash J(A,d)\lra U(\Lambda)\backslash J(A)$.
Again, since $(A^\times\cap\ker(d))\subseteq A^\times$ the same map yields a
group homomorphism
$$U(\Lambda,d)\backslash J(A,d)/(A^\times\cap\ker(d))\lra U(\Lambda)\backslash J(A)/A^\times$$
and since $Cl(\Lambda)\simeq U(\Lambda)\backslash J(A)/A^\times$, we get the statement. \dickebox

\begin{Example} \label{notsurjectivemapofclassgroups}
Consider the $\Z$-order
$$\Lambda=\Z\cdot\left(\begin{array}{cc}1&0\\0&1\end{array}\right)+p\cdot Mat_{2\times 2}(\Z)$$
in $Mat_{2\times 2}(\Q)$. Then, following Example~\ref{Examplediffgradedalgebra}
this is a differential graded order with the differential $d_1$.
{\red Indeed, for $a,x,y,z,u\in\Z$ we compute
$$d_1(\left(\begin{array}{cc}a+px&pu\\pz&a+py\end{array}\right)=
\left(\begin{array}{cc}pz&p(y-x)\\0&pz\end{array}\right)\in\Lambda$$}
Further, by \cite[\S 34D]{CR1} we have that $Cl(\Lambda)$ is of order $2$
whenever $p-1$ is divisible by $4$ (and trivial otherwise).
However, $\ker(d_1)=\Z[\epsilon]/\epsilon^2$, and hence $Cl(\Lambda,d_1)$, {\red as
well as the dg-id\`ele group are} trivial
as is easily seen using the id\`ele description. Therefore, the map
in Proposition~\ref{ordinaryclassgroupdgclassgroup} is not surjective in general.
\end{Example}

Since localisation is flat, we get $(H(\Lambda,d))_\wp\simeq H(\Lambda_\wp,d_\wp)$. For any
dg-module $(M,\delta)$ we get $(H(M,\delta))_\wp\simeq H(M_\wp,\delta_\wp)$. Moreover,
if $(M,\delta)$ is a differential graded $(\Lambda,d)$-module, then by  Lemma~\ref{homologymodule}
we get that  $H(M,\delta)$ is a graded $H(\Lambda,d)$-module.
Further, if $(M,\delta)$ is a locally free dg $(\Lambda,d)$-module, then
$(M_\wp,\delta_\wp)$ is a free $(\Lambda_\wp,d_\wp)$-module.
Hence $H(M,\delta)$ is a graded locally free $H(\Lambda,d)$-module.

If we identify $$Cl(H(\Lambda,d))=\ker\left({\red G}_0(H(\Lambda,d)-LF-mod)\lra {\red G}_0(H(A,d)-proj)\right)$$
taking homology is therefore a group homomorphism
$$Cl(\Lambda,d)\stackrel{CH_{(\Lambda,d)}}\lra Cl(H(\Lambda,d)).$$

Note however that $H(\Lambda,d)$ is not an order in general.

\begin{Def}
Let $R$ be the ring of integers in a number field $K$. Let $(\Lambda,d)$ be a
differential graded order in the
semisimple differential graded $K$-algebra $(A,d)$.
Then define the {\em homology-isomorphism class group kernel} $\ker(CH_{(\Lambda,d)})=:Cl(\Lambda,d)_{hi}$.
\end{Def}

\begin{Rem}
Note that for each locally free differential graded $(\Lambda,d)$-module $(L,\delta)$ we have that
$$Cl(\Lambda,d)_{hi}+[(L,\delta)]=CH_{(\Lambda,d)}^{-1}(CH_{(\Lambda,d)}([(L,\delta)])$$
parameterizes those locally free dg-$(\Lambda,d)$-modules which have the same homology
as $(L,\delta)$. Those which are quasi-isomorphic to $(L,\delta)$ do share the homology with
$(L,\delta)$. However, being quasi-isomorphic is stronger in general than just having isomorphic
homology.
\end{Rem}

We now consider locally free $H(\Lambda,d)$-modules.

\begin{Lemma}\label{torsionsubmodulesoflocallfree}
Let $R$ be a Dedekind domain with field of fractions $K$ and let
$B$ be a finitely generated $R$-algebra. For any $B$-module $M$ let
$$t(M):=\{x\in M\;|\;\exists r\in R\setminus\{0\}:rx=0\}$$
be the torsion submodule of $M$. If $L$ is a locally free $B$-module of rank $n$,
then $t(L)\simeq t(B)^n$.
\end{Lemma}

Proof. By hypothesis $L_\wp\simeq B_\wp^n$ for all primes $\wp\in Spec(R)$.
Since $R$ is Dedekind, $R_\wp$ is a principal ideal domain (cf \cite[Lemma 7.5.9]{grouprepbuch}).
Hence $t(L_\wp)=t(L)_\wp$. Since $R$ is Dedekind, any non zero prime ideal is maximal
(cf \cite[Lemma 7.5.15]{grouprepbuch}) and hence the primary
decomposition (cf \cite[Theorem 7.2.5]{grouprepbuch})
of $t(L)$ gives a direct product decomposition
$t(L)=\prod_{\wp\in Spec(R)}t(L)_\wp$.
We hence get $t(L)\simeq t(B)^n$. \dickebox

\begin{Cor}
Let $R$ be a Dedekind domain with field of fractions $K$, let $(A,d)$ be a
finite dimensional and semisimple and differential graded
$K$-algebra. Let $(\Lambda,d)$ be a differential graded $R$-order
in $(A,d)$.
Then $$Cl(H(\Lambda,d))\simeq Cl(H(\Lambda,d)/t(H(\Lambda,d))).$$
In particular, if $H(A,d)=0$, then  $Cl(\Lambda,d)_{hi}=0$.
\end{Cor}

Proof. By Lemma~\ref{torsionsubmodulesoflocallfree} for any rank $1$ locally free
$H(\Lambda,d)$-module $L$ we have $t(L)=t(H(\Lambda,d))$ and therefore
$$H(\Lambda,d)/t(H(\Lambda,d))\otimes_{H(\Lambda,d)}-:Cl(H(\Lambda,d))\lra Cl(H(\Lambda,d)/t(H(\Lambda,d)))$$
is an isomorphism. Note that $H(\Lambda,d)/t(H(\Lambda,d))\otimes_{H(\Lambda,d)}-\simeq (R/tR)\otimes_R-$.

If now $(A,d)$ is {\red acyclic}, then $H(\Lambda,d)$ is torsion, and
therefore locally free $H(\Lambda,d)$-modules are actually free.
This proves the lemma. \dickebox

\medskip

We shall need a result due to Wehlen. Recall (cf e.g. \cite[Chapter 3]{GoodearlWarfield})
that the Baer lower radical $L(A)$ of an algebra $A$
is the intersection of all prime ideals of $A$. Every element of $L(A)$ is nilpotent
and if $A$ is Noetherian, then $L(A)$ is a nilpotent ideal \cite[Theorem 3.11]{GoodearlWarfield}.

\begin{Theorem} \cite[Theorem 2.4]{Wehlen} \label{Wehlentheorem}
Let $R$ be a Pr\"ufer domain and let $A$ be a finitely generated algebra over $R$.
Let $L(A)$ be the (Baer) lower radical of $A$. If $B:=A/L(A)$ is separable over $R$,
i.e. $B$ is projective as an $B\otimes_RB^{op}$-module, then there are idempotents
$e_p$ and $e_t$ of $A$ such that
$$A=\left(\begin{array}{cc}e_pAe_p&e_tAe_p\\e_pAe_t&e_tAe_t\end{array}\right)$$
and $e_pAe_p/e_pL(A)e_p$ is $R$-projective, such that $\left(\begin{array}{cc}0&e_tAe_p\\e_pAe_t&e_tAe_t\end{array}\right)\subseteq t(A)$,
the $R$-torsion ideal of $A$.
\end{Theorem}

Note that obviously $e_p$ and $e_t$ are orthogonal and $e_p+e_t=1$. Further, note that
$A/t(A)=e_pAe_p/t(e_pAe_p)$.

\begin{Cor} \label{liftingunitsalaWehlen}
Suppose 
the hypotheses of Theorem~\ref{Wehlentheorem}. Then any unit $u$ in $A/t(A)$
can be lifted to a unit $u_2$ of $A$.
\end{Cor}

Proof.
If $u$ is a unit of $A/t(A)$, then we get that $u$ is actually
a unit of $e_pAe_p/t(e_pAe_p)$. Further, since $e_pAe_p/e_pL(A)e_p$ is $R$-projective,
we have that $t(e_pAe_p)\subseteq e_pL(A)e_p$ and since the right hand side
of the inclusion is a nil ideal, so is the left hand side.
Hence there is $u_1$ and $v_1$ in $e_pAe_p$ such that $u_1$ maps to $u$ in $A/t(A)$ and such that
$u_1v_1-1$ is in $e_pL(A)e_p$, whence nilpotent. But since $1+n$ for a nilpotent element $n$ is a
unit, we have that $u_1$ is a unit. Therefore we may find a unit
$u_1\in e_pAe_p$ such that $u_1$ maps to $u$ in $A/t(A)$. But then there is a unit
$u_2:=\left(\begin{array}{cc}u_1&0\\ 0&1\end{array}\right)$ in $A$ (actually in the Pierce
decomposition above) such that $u_2$ maps to $u$ in $A/t(A)$. \dickebox

\begin{Theorem}
Let $R$ be a Dedekind domain with field of fractions $K$ and let $(A,d)$ be a differential graded
algebra such  that $A$ is finite dimensional separable over $K$. Suppose that $\ker(d)$ is separable as
a graded algebra.
Let $(\Lambda,d)$ be a differential graded
$R$-order in $(A,d)$. Denote the canonical map $\pi:\ker(d)\lra H(\Lambda,d)$ and denote $\overline{H(\Lambda,d)}:=H(\Lambda,d)/t(H(\Lambda,d))$.
Then $\overline{H(\Lambda,d)}$ is a classical
$R$-order in $H(A,d)$, and we have an exact sequence
$$0\lra Cl(\Lambda,d)_{hi}\lra Cl(\Lambda,d)\lra Cl(\overline{H(\Lambda,d)})\lra  0.$$
\end{Theorem}

Proof. It is a classical fact that separable algebras are semisimple.
Since $\ker(d)$ is a semisimple algebra, Lemma~\ref{semisimplekerd} shows that
we see that $H(A,d)$ is a semisimple
algebra. Hence, $H(\Lambda,d)/t(H(\Lambda,d))$ is an $R$-order in the semisimple
$K$-algebra $H(A,d)$.

We need to see that the right hand map is surjective. For this we shall use the interpretation
of class groups as id\`eles.

We hence need to see that
$$\widehat U(\Lambda,d)\backslash \widehat J(A,d)/(A^\times\cap\ker(d))\lra
\widehat U(\overline{H(\Lambda,d)})\backslash \widehat J(H(A,d))/H(A,d)^\times$$
is surjective, {\red since the map on the level of id\`eles factors through the right hand map}.

First, by Lemma~\ref{semisimplekerd} any unit of $H(A,d)$ lifts to
a unit in $\ker(d)$. Now, clearly $\ker(d)^\times\subseteq(A^\times\cap\ker(d))$.
For a finitely generated $R$-algebra $B$ satisfying that $B/L(B)$ is separable,
if $\ol u$ is a unit in $B/tB$, then by Corollary~\ref{liftingunitsalaWehlen}
any unit $u$ in $B/tB$ lifts to a unit $u$ in $B$. Now, by hypothesis $\ker(d)$ is
finite dimensional separable. This is equivalent to the fact that
in the Wedderburn decomposition the centres of the skew fields are separable
extensions of the base field.
Hence $H(A,d)$ is separable as well. Therefore $H(\Lambda_\wp,d)/t(H(\Lambda_\wp,d))$
is separable for all $\wp\in Spec(R)$, and so is its quotient modulo its Baer (lower) radical.
Following Corollary~\ref{liftingunitsalaWehlen} we can hence lift any unit of
$H(\Lambda_\wp,d)/t(H(\Lambda_\wp,d))$ to a unit of
$H(\Lambda_\wp,d)$


By semisimplicity $A_\wp^\times\cap\ker(d)\lra H(A_\wp,d)^\times$ is surjective.
By the same argument $\widehat J(A,d)\lra \widehat J(H(A,d))$ is surjective as well.
Since $\pi(\widehat U(\Lambda,d))\subseteq \widehat U(H(\Lambda,d))$, we get indeed a  map
$Cl(\Lambda,d)\lra Cl(\overline{H(\Lambda,d)})$.

Now, for any locally free $\ol{H(\Lambda,d)}$-ideal $\ol I$ we can find an id\`ele
$\ol \alpha\in J(\ol{H(A,d)})$ such that $\ol I=\bigcap_{\wp\in Spec (R)}\ol{H(\Lambda,d)}\cdot\ol \alpha_\wp$.
Since $\widehat J(A,d)\lra \widehat J(H(A,d))$ is surjective, there is an id\`ele $\alpha\in J(A,d)$
which maps to $\ol\alpha$. But then
the locally free ideal
$$I:=\bigcap_{\wp\in Spec (R)}\Lambda\cdot\alpha_\wp$$
maps to $\ol I$. This shows the surjectivity of $Cl(\Lambda,d)\lra Cl(\overline{H(\Lambda,d)})$.
The statement on the kernel of this map follows by definition.
\dickebox


%

%
%
%
%

\section{Reducing to dg-orders in differential graded simple algebras}

\label{mayervietorissection}

A major tool in studying class groups of orders is the following result.

\begin{Theorem} (Reiner-Ullom~\cite{RU74})(Mayer-Vietoris theorem for class
groups of orders)
\label{classgrouppullback}
Let $R$ be a Dedekind domain with field of fractions $K$ and let $\Lambda$ be an
$R$-order in the finite dimensional semisimple $K$-algebra $A$ satisfying the Eichler condition.
Let $e^2=e\in Z(A)$ be a
non trivial central idempotent of $A$. Put $f:=1-e$. Then $\Lambda e$ is an $R$-order
in $Ae$ and $\Lambda f$ is an $R$-order in $Af$ and we have a pullback diagram
$$
\xymatrix{\Lambda\ar[r]^{\cdot e}\ar[d]_{\cdot f}&\Lambda\cdot e\ar[d]^{\pi_e}\\
\Lambda f\ar[r]_{\pi_f}&\ol{\Lambda}
}
$$
for $\ol\Lambda:=\Lambda e/\Lambda\cap\Lambda e\simeq \Lambda f/\Lambda\cap\Lambda f$
and $\pi_e$, resp. $\pi_f$ being the canonical morphisms.
Further, there is a group homomorphism $\delta$ such that $\delta$ and the canonical maps induce
an exact sequence
$$
\Lambda^\times\lra(\Lambda e)^\times\times(\Lambda f)^\times\lra\ol\Lambda^\times\stackrel\delta\lra
Cl(\Lambda)\lra Cl(\Lambda e)\times Cl(\Lambda f).$$
\end{Theorem}

Keep the assumptions and notations of Theorem~\ref{classgrouppullback} for the moment,{ \red and suppose that $K$ is of characteristic different from $2$.} In addition
if $(A,d)$ is a finite dimensional semisimple {\red (as an algebra)} differential graded $K$-algebra, and
if $(\Lambda,d)$ is a differential graded $R$-order in $(A,d)$, then
by Remark~\ref{blocksofdgalgebras} we have $d(e)=0=d(f)$
and hence $(\Lambda e,d)$ and $(\Lambda f,d)$ are differential graded orders in
$(Ae,d)$ respectively $(Af,d)$.
Since $\ol\Lambda=\Lambda e/(\Lambda\cap\Lambda e)$, and since $d$ induces a differential on both
$\Lambda$ and on $\Lambda e$, it also induces a differential $\ol d$ on $\ol\Lambda$.
Therefore, the maps in
$$
\xymatrix{(\Lambda,d)\ar[r]^{\cdot e}\ar[d]_{\cdot f}&(\Lambda\cdot e,d)\ar[d]^{\pi_e}\\
(\Lambda \cdot f,d)\ar[r]_{\pi_f}&(\ol{\Lambda},\ol d)
}
$$
are still well-defined maps of differential graded algebras. Since
$$
\xymatrix{\Lambda\ar[r]^{\cdot e}\ar[d]_{\cdot f}&\Lambda\cdot e\ar[d]^{\pi_e}\\
\Lambda \cdot f\ar[r]_{\pi_f}&\ol{\Lambda}
}
$$
is a pullback diagram of $R$-algebras,
the diagram
$$
\xymatrix{(\Lambda,d)\ar[r]^{\cdot e}\ar[d]_{\cdot f}&(\Lambda\cdot e,d)\ar[d]^{\pi_e}\\
(\Lambda \cdot f,d)\ar[r]_{\pi_f}&(\ol{\Lambda},\ol d)
}
$$
is a pullback diagram of differential graded algebras.
Recall that $Cl(\Lambda,d)$ is generated by $[L]-[\Lambda]$ for $L$ being
rank $1$ locally free differential graded modules. Hence the map
$$ Cl(\Lambda,d)\stackrel{Cl(e)\times Cl(f)}\lra Cl(\Lambda e,d)\times Cl(\Lambda f,d)$$
is still well-defined. We recall the definition of $\delta$.
Let $u\in\ol\Lambda^\times$. Then consider the pullback diagram
$$
\xymatrix{L_u\ar[rr]^{\alpha_e}\ar[dd]_{\alpha_f}&&\Lambda\cdot e\ar[d]^{\pi_e}\\
&&\ol{\Lambda}\ar[dl]^{\cdot u}\\
\Lambda \cdot f\ar[r]_{\pi_f}&\ol{\Lambda}}
$$
and define $\delta(u):=[L_u]-[\Lambda]$.

In the dg-case, if $u\in\ker(\ol d)\cap\ol\Lambda^\times$  {\red is homogeneous (and hence of degree $0$)}, then we get again a pullback diagram
$$
\xymatrix{(L_u,d_u)\ar[rr]^{\alpha_e}\ar[dd]_{\alpha_f}&&(\Lambda\cdot e,d)\ar[d]^{\pi_e}\\
&&(\ol{\Lambda},\ol d)\ar[dl]^{\cdot u}\\
(\Lambda \cdot f,d)\ar[r]_{\pi_f}&(\ol{\Lambda},d)}
$$
defining a locally free differential graded $(\Lambda,d)$-module $(L,d_u)$.
We may define a map $\delta_d$ by
$\delta_d(u):=[(L_u,d_L)]-[(\Lambda,d)]$ in this case.
It remains to show that $\ker(Cl(\Lambda e)\times Cl(\Lambda f))=\im(\delta_d)$.

Suppose that $u\in\ker(\ol d)\cap\ol{\Lambda}^\times$. Then
$$Cl(f)\circ\delta_d(u)=[(L_u\cdot f,d_u)]-[(\Lambda f,d)]=[\alpha_f(L_u,d_u)]-[(\Lambda f,d)]=
[(\Lambda f,d)]-[(\Lambda f,d)]=0
$$
and likewise $Cl(e)\circ\delta_d(u)=0$.

Let $(L,d_L)$ be a rank $1$ locally free dg $(\Lambda,d)$-module in the kernel of $Cl(e)\times Cl(f)$.
Let $\alpha$ be the id\`ele of $(L,d_L)$. Then we know that
$\alpha\cdot e$ is the principal id\`ele and also $\alpha\cdot f$ is the principal id\`ele
in the corresponding algebras. Now
$$A^\times\cap\ker d=(Ae^\times\cap\ker d)\times (Af^\times\cap\ker d).$$
Further, $\ker(d)=\ker(de)\times\ker(df)$. We hence only need to  consider ordinary id\`eles
in order to examine the kernel of $Cl(e)\times Cl(f)$. But then, by the classical Mayer-Vietoris-like
theorem for class groups of orders we obtain that
$$\ker(Cl(e)\times Cl(f))=\im(\delta_d).$$

We hence obtained the following

\begin{Theorem}\label{MayerVietorisfordg}
Let $R$ be a Dedekind domain with field of fractions $K$
{\red of characteristic different from $2$} and let $(\Lambda,d)$ be a
differential graded
$R$-order in the finite dimensional semisimple {\red (as an algebra) }
differential graded $K$-algebra $(A,d)$,
and suppose that $A$ satisfies the Eichler condition.
Let $e^2=e\in Z(A)$ be a
non trivial central idempotent of $A$. Put $f:=1-e$. Then $(\Lambda e,de)$
is a differential graded $R$-order
in $(Ae,de)$ and $(\Lambda f,df)$ is a differential graded $R$-order in $(Af,df)$
and we have a pullback diagram
$$
\xymatrix{\Lambda\ar[r]^{\cdot e}\ar[d]_{\cdot f}&\Lambda\cdot e\ar[d]^{\pi_e}\\
\Lambda\cdot f\ar[r]_{\pi_f}&\ol{\Lambda}_e
}
$$
for $\ol\Lambda_e:=\Lambda e/\Lambda\cap\Lambda e\simeq \Lambda f/\Lambda\cap\Lambda f$
and $\pi_e$, resp. $\pi_f$ being the canonical morphisms.
Further, there is a group homomorphism $\delta_d$ such that $\delta_d$
and the canonical maps induce an exact sequence
$$
\Lambda^\times\cap\ker(d)\lra((\Lambda e)^\times\times(\Lambda f)^\times)\cap\ker(d)\lra\ol\Lambda_e^\times\cap\ker(\ol d)\stackrel{\delta_d}\lra
Cl(\Lambda,d)\lra Cl(\Lambda e,de)\times Cl(\Lambda f,df).$$
\end{Theorem}

Proof. The result follows from the comments leading to the
exactness of
$$\ol\Lambda_e^\times\cap\ker(\ol d)\stackrel{\delta_d}\lra
Cl(\Lambda,d)\lra Cl(\Lambda e,de)\times Cl(\Lambda f,df).$$
Since $\Lambda$ is an order in a semisimple algebra $A$ and $e^2=e\in Z(A)$,
the sequence
$$
\Lambda^\times\stackrel{\gamma}\lra((\Lambda e)^\times\times(\Lambda f)^\times)
\stackrel{\pi}\lra\ol\Lambda_e^\times
$$
is exact by the classical situation.
Since
$$
\xymatrix{\Lambda\ar[r]^{\cdot e}\ar[d]_{\cdot f}&\Lambda\cdot e\ar[d]^{\pi_e}\\
\Lambda \cdot f\ar[r]_{\pi_f}&\ol{\Lambda}_e
}
$$
is a pullback diagram, we get an induced sequence of group homomorphisms
$$
\Lambda^\times\cap\ker(d)\lra((\Lambda e)^\times\times(\Lambda f)^\times)\cap\ker(d)\stackrel{\pi_d}\lra\ol\Lambda_e^\times\cap\ker(\ol d)
$$
which composes to $0$, and where $\pi_d$ is the restriction of $\pi$ to
$((\Lambda e)^\times\times(\Lambda f)^\times)\cap\ker(d)$.
If $x\in\ker(\pi_d)$ {\red is homogeneous}, then there is {\red a homogeneous}
$v\in\Lambda^\times$ with $\gamma(v)=x$
{\red (where $\gamma:\Lambda^\times\lra (\Lambda e)^\times\times(\Lambda f)^\times$).}
We need to see that $v\in\ker(d)$.
However, $K\Lambda=K\Lambda e\times K\Lambda f$ and $K\ol\Lambda=0$.
Hence, $$\textup{id}_K\otimes_R \gamma=\textup{id}_A$$
which shows that $x\in\ker(d)$ implies $v\in\ker(d)$.
This proves  the theorem. \dickebox

\begin{Rem}\label{comparewithhomologyofodLambda}
Taking homology we can compare the Mayer-Vietoris sequence of Theorem~\ref{MayerVietorisfordg}
with the Mayer-Vietoris sequence of Theorem~\ref{classgrouppullback}. Denoting as before
$$\ol{H(\Lambda,d)}=H(\Lambda,d)/tH(\Lambda,d),$$
and likewise for the other orders occurring in the statement,
observing that there is a
map $\ker(\ol d)\lra H(\ol\Lambda_e,\ol d)$, we get a commutative diagram
$$
\xymatrix{((\Lambda e)^\times\times(\Lambda f)^\times)\cap\ker(d)\ar[r]\ar[d]
&\ol\Lambda_e^\times\cap\ker(\ol d)\ar[r]^{\delta_d}\ar[d]& Cl(\Lambda,d)\ar[r]\ar[d]&
Cl(\Lambda e,de)\times Cl(\Lambda f,df)\ar[d]\\
\ol{H(\Lambda e,de)}^\times\times \ol{H(\Lambda f,df)}^\times\ar[r]&
\ol{H(\ol\Lambda_e,\ol d)}^\times\ar[r]^{\ol\delta}&Cl(\ol{H(\Lambda,d)})\ar[r]&
Cl(\ol{H(\Lambda e,de)})\times Cl(\ol{H(\Lambda f,df)})
}
$$
\end{Rem}

In the situation of Theorem~\ref{MayerVietorisfordg} we shall need a particular subgroup
of the unit group of $\ol\Lambda$.

\begin{Def}
Let $R$ be a Dedekind domain with field of fractions $K$ {\red of characteristic different from $2$}
and let $(\Lambda,d)$ be a
differential graded
$R$-order in the finite dimensional semisimple differential graded $K$-algebra $(A,d)$.
Let $e^2=e\in Z(A)$ and $f=1-e$.
Let $D_e$ be the image of
$((\Lambda e)^\times\times(\Lambda f)^\times)\cap\ker(d)\lra\ol\Lambda^\times\cap\ker(\ol d)$.
Put
$$X_e:=\ker((\ol\Lambda_e^\times\cap\ker \ol d)\ra\ol{ H(\ol \Lambda_e,\ol d)}^\times)$$
we define homology-isomorphism units of $\ol\Lambda_e:=\Lambda e/(\Lambda\cap\Lambda e)$ as
$$\ol\Lambda_{hi}^\times(e):=
X_e/\langle [X_e,X_e], D_e\rangle.$$
\end{Def}

\begin{Prop}\label{shortexactsequenceforqis}
Let $R$ be a Dedekind domain with field of fractions $K$ {\red of characteristic different from $2$}
and let $(\Lambda,d)$ be a
differential graded
$R$-order in the finite dimensional semisimple differential graded $K$-algebra $(A,d)$.
Suppose that $H(A,d)$ is semisimple as an algebra.
{\red Suppose that $A$ and $H(A,d)$ both satisfy the Eichler condition.}
Let $e^2=e\in Z(A)$ be a
non trivial central idempotent of $A$. Put $f:=1-e$ and $\ol\Lambda:=\Lambda e/(\Lambda\cap\Lambda e)$.
Then the sequence
$$
\ol\Lambda^\times_{hi}(e)\lra Cl_{hi}(\Lambda,d)\lra Cl_{hi}(\Lambda e,de)\times Cl_{hi}(\Lambda f,df)
$$
is exact.
\end{Prop}

Proof. The proof only uses standard diagram chasing. Here are the details.
Recall from Remark~\ref{comparewithhomologyofodLambda} that there is a commutative diagram
with exact rows
$$
\xymatrix{((\Lambda e)^\times\times(\Lambda f)^\times)\cap\ker(d)\ar[r]\ar[d]
&\ol\Lambda_e^\times\cap\ker(\ol d)\ar[r]^{\delta_d}\ar[d]& Cl(\Lambda,d)\ar[r]\ar[d]&
Cl(\Lambda e,de)\times Cl(\Lambda f,df)\ar[d]\\
\ol{H(\Lambda e,de)}^\times\times \ol{H(\Lambda f,df)}^\times\ar[r]&\ol{H(\ol\Lambda,\ol d)}^\times\ar[r]^{\ol\delta}&Cl(\ol{H(\Lambda,d)})\ar[r]&
Cl(\ol{H(\Lambda e,de)})\times Cl(\ol{H(\Lambda f,df)}).
}
$$
Since the class groups are abelian groups, we can consider the abelianisations of the various unit groups,
denoted by the index $ab$, for short.
Put $F_e$ the image of
$$\ol{H(\Lambda e,de)}^\times\times \ol{H(\Lambda f,df)}^\times\lra\ol{H(\ol\Lambda,\ol d)}^\times.$$
This then gives a commutative diagram with exact rows (where we replace $D_e$ and $F_e$ by the image in the abelianisation
$$
\xymatrix{
\left(\ol\Lambda_e^\times\cap\ker(\ol d)\right)_{ab}/D_e\ar@{^{(}->}[r]^{\delta_d}\ar[d]& Cl(\Lambda,d)\ar[r]\ar[d]&
Cl(\Lambda e,de)\times Cl(\Lambda f,df)\ar[d]\\
\left(\ol{H(\ol\Lambda_e,\ol d)}^\times\right)_{ab}/F_e\ar@{^{(}->}[r]^{\ol\delta}&Cl(\ol{H(\Lambda,d)})\ar[r]&
Cl(\ol{H(\Lambda e,de)})\times Cl(\ol{H(\Lambda f,df)})
}$$
Taking kernels of the vertical maps can be extended to a commutative diagram with exact columns
$$
\xymatrix{
\ol\Lambda^\times_{hi}\ar[r]^-{\alpha_K}\ar@{^{(}->}[d]^{\iota_1}&
Cl_{hi}(\Lambda,d)\ar[r]^-{\beta_K}\ar@{^{(}->}[d]^{\iota_2}&
Cl_{hi}(\Lambda e,de)\times Cl_{hi}(\Lambda f,df)\ar@{^{(}->}[d]^{\iota_3}\\
\left(\ol\Lambda_e^\times\cap\ker(\ol d)\right)_{ab}/D_e\ar@{^{(}->}[r]^-{\delta_d}\ar[d]^{\pi_1}& Cl(\Lambda,d)\ar[r]^-{\epsilon}\ar[d]^{\pi_2}&
Cl(\Lambda e,de)\times Cl(\Lambda f,df)\ar[d]^{\pi_3}\\
\left(\ol{H(\ol\Lambda_e,\ol d)}^\times\right)_{ab}/F_e\ar@{^{(}->}[r]^-{\ol\delta}&Cl(\ol{H(\Lambda,d)})\ar[r]&
Cl(\ol{H(\Lambda e,de)})\times Cl(\ol{H(\Lambda f,df)})
}
$$
Clearly $\beta_K\circ\alpha_K=0$, which shows $\im(\alpha_K)\subseteq\ker(\beta_K)$.

If $x\in\ker(\beta_K)$, then
$$\iota_2(x)=\delta_d(y)\in\ker(\epsilon)\cap\ker(\pi_2)$$
for some $y\in \ol\Lambda^\times\cap\ker(\ol d)$ and
$$(\ol\delta\circ\pi_1)(y)=(\pi_2\circ\delta_d)(y)=\pi_2(x)=0.$$
Since $\ol\delta$ is injective, $\pi_1(y)=0$, which shows that there is $y_1\in\ol\Lambda^\times_{hi}$
such that $\iota_1(y_1)=y$. Then
$$(\iota_2\circ\alpha_K)(y_1)=\delta_d(\iota_1(y_1))=\delta_d(y)=\iota_2(x).$$
Since $\iota_2$ is injective, $\alpha_K(y_1)=x$. This shows the statement. \dickebox

\subsection*{Acknowledgements} I wish to thank Bernhard Keller for encouragement in an
early stage of the project, and for useful comments on the first draft of the manuscript.
Thanks to Yann Palu for indicating reference \cite{Raedschelder} to me.
I am very grateful to the referee for extremely
careful reading. The referee spotted numerous imprecisions, omissions,
and mistakes, in particular in Section~\ref{ringtheoryofdgalgebras}. This lead to very substantial
improvements of the paper.

\end{document}